\documentclass[11pt,a4paper]{article}
\usepackage[latin1]{inputenc}
\usepackage[T1]{fontenc}
\usepackage[francais]{babel}
\usepackage{makeidx}
\usepackage{amsfonts}
\usepackage{amsthm}
\usepackage{amssymb}
\usepackage{amsmath}
\usepackage[all]{xy}
\begin{document}

\begin{center}
\Large \textbf{Nombres de Weil, sommes de Gauss et annulateurs Galoisiens}
\end{center}
\begin{center}
\textbf{Thong NGUYEN QUANG DO, Vésale NICOLAS}
\end{center}

\tableofcontents
\ \\
\textbf{Abstract:}

For an abelian number field $K$ containing a primitive $p^{th}$ root of unity ($p$ an odd prime) and satisfying certain
 technical conditions, we parametrize the $\mathbb{Z}_p[\mathrm{G}(K/\mathbb{Q})]$-annihilators of the "minus" part $A_{K}^{-}$ of the $p$-class group by means of modules of Jacobi sums. Using a reflection theorem and Bloch-Kato's reciprocity law, we then determine the Fitting ideal of the "plus" part $A_{K}^{+}$ in terms of "twisted" Gauss sums.
\\ \ \\
\textbf{Key words:} sommes de Gauss, miroir de Leopoldt, représentation d'Ihara, réciprocité de Bloch-Kato.
\\ \ \\
\textbf{AMS subject classification:} 11R23.

\section{Introduction}

Pour une extension galoisienne de corps de nombres $K/k$, de groupe de Galois $H$, l'étude de la structure de $H$-module du groupe de classes $\mathcal{C}l_{K}$ commence naturellement par la détermination de son idéal annulateur. Si $k=\mathbb{Q}$ et $K$ est un corps abélien CM, un annulateur classique est l'idéal de Stickelberger, mais qui ne voit que "la moitié du monde", en ce sens que sa partie "plus" est triviale, et donc n'apporte que des renseignements... triviaux. Les mêmes problèmes se posent dans le cadre $p$-adique pour la $\mathbb{Z}_p[H]$-structure du $p$-groupe de classes $A_K$: la partie "moins" s'étudie à partir de la Conjecture Principale de la théorie d'Iwasawa (théorème de Mazur-Wiles), mais la partie "plus" demeure mystérieuse (conjecture de Vandiver, conjecture de Greenberg).
\\Dans le présent travail, on fixe un corps abélien $K$ contenant une racine primitive $p^{ieme}$ de l'unité $\zeta_p$ ($p\neq 2$) et l'on adopte les notations habituelles $K_n=K(\zeta_{p^{n+1}})$, $K_{\infty}=\cup_{n\geq 0} K_n$, $H_n=\mathrm{Gal}(K_n/\mathbb{Q})$, $H_{\infty}=\mathrm{Gal}(K_{\infty}/\mathbb{Q})$, $A_n$=le $p$-groupe de classes de $K_n$, $X_{\infty}=\displaystyle{\lim_{\longleftarrow}}\ A_n$. On se propose dans un premier temps de "paramétriser" les $\mathbb{Z}_p[H_n]$ (resp. $\mathbb{Z}_p[[H_{\infty}]]$) annulateurs de $A_n^{-}$ (resp. $X_{\infty}^{-}$) par des modules de sommes de Jacobi: c'est une idée ancienne d'A. Weil ([We]), appliquée par Iwasawa à $\mathbb{Q}(\zeta_p)$ ([Iw1]), puis par Anglès et Beliaeva aux corps $\mathbb{Q}(\zeta_{p^n})$ et $\mathbb{Q}(\zeta_{p^{\infty}})$ ([AB]), et généralisée aux corps abéliens par le second auteur dans son mémoire de Master. Dans un second temps, on utilise le "miroir" (\textsl{Spiegelung}) pour en déduire une description de l'idéal de Fitting initial de $A_n^{+}$ (resp. $X_{\infty}^{+}$) sur $\mathbb{Z}_p[H_n]$ (resp. $\mathbb{Z}_p[[H_{\infty}]]$) en termes de sommes de Gauss "tordues" (théorème 3.4.2.). Ces sommes de Gauss sont explicitées dans une troisième étape en utilisant la loi de réciprocité de Bloch-Kato (théorème 5.3.2. ; \textit{le principe du calcul est expliqué à la fin du §3}). Dans un dernier paragraphe, on fait le lien avec des résultats d'annulation récemment annoncés par D.Solomon ([Sol1]) au colloque "Iwasawa 2008" à Irsee. Un appendice rappelle les relations du "miroir" (dont certaines sont peut-être mal connues) qui sont utilisées dans le texte.
\\ \ \\Dans la majeure partie de l'article, les calculs sont faits caractère par caractère, et l'on se place pour simplifier dans le cas où $p$ ne divise pas $[K:\mathbb{Q}]$. Ce n'est pas une contrainte essentielle, le cas général pouvant se traiter par exemple par les méthodes de [T1] et [B], la plupart des résultats demeurant valables, mais au prix de techniques plus compliquées et d'énoncés moins précis (\textit{voir les commentaires du théorème 3.3.1. du texte}). Par ailleurs, puisque $p\neq 2$, on pourrait aussi se concentrer sur le caractère quadratique de $\mathrm{Gal}(K/K^{+})$ (sans hypothèse de semi-simplicité) pour obtenir, avec les méthodes du présent travail, des résultats "globaux" sur les parties "+ et -" des modules et idéaux concernés. Mais qu'on se place dans la situation semi-simple ou pas, les caractères impairs $\psi$ tel que $\psi(p)=1$ posent problème: c'est le cas dit "totalement $p$-décomposé", notoirement réfractaire aux techniques de descente de la théorie d'Iwasawa, et où les complications ne sont plus seulement techniques. Un exemple typique est le théorème 3.3.1. ci-après, où le résultat est général au niveau infini, mais ne "descend bien" au niveau fini que sous une hypothèse de non-décomposition. On a écarté le cas $p$-décomposé de la descente pour ne pas surcharger cet article (voir les commentaires du §3.4).

\ \\
\textbf{Remerciements}
\\ \ \\
Nos remerciements vont à Anthony Martin pour sa lecture attentive du manuscrit et des discussions utiles sur la partie technique des lois de réciprocité, ainsi qu'au rapporteur dont les remarques et suggestions nous ont permis d'améloirer la présentation et l'exposition de ce travail.

\section{Nombres de Weil et sommes de Gauss}
\noindent{\textbf{Notations générales}}
\\ \ \\
Soit $p$ un nombre premier impair. 
Soit $K/\mathbb{Q}$ une extension abélienne de $\mathbb{Q}$ contenant $\mu_p$. Soit $F$ un sous-corps totalement réel de $K$. 
\\Notons, pour $n\geq1$ $K_n=K(\mu_{p^{n+1}})$, $G_n=\mathrm{Gal}(K_n/F)$, $\Gamma_n=\mathrm{Gal}(K_n/K)$, $\Delta=\mathrm{Gal}(K/F)$. 
\\Notons $\mu_{K_n}$ le groupe des racines de l'unité de $K_n$ et $c_n$ le conducteur de $K_n$.
\\Soit $I_n$ le $\mathbb{Z}[G_n]$-module des idéaux fractionnaires de $K_n$ premiers à $c_0$, $P_n$ le $\mathbb{Z}[G_n]$-module des $\alpha\in K_n$ tels que $(\alpha)\in I_n$,  $\mathcal{O}_{K_n}$ l'anneau des entiers de $K_n$ et $\mathcal{C}l_{K_n}$ le groupe des classes de $K_n$.
Notons $\mathcal{U}_n$ le module des unités semi-locales de $K_n$, i.e.
$$\mathcal{U}_n=\displaystyle{\prod_{\mathfrak{q}|p}}U_{\mathfrak{q},n}^{1};$$
où $U_{\mathfrak{q},n}^{1}$ est le module des unités principales du $\mathfrak{q}$-complété de $K_n$.

\ \\

Dans ce paragraphe, nous généralisons les définitions des modules et des nombres de Weil de la première partie de [AB] au cas abélien sur $\mathbb{Q}$ et nous relions ces modules à l'annulateur du groupe de classes de $K_n$, aux sommes de Gauss ([Ich1],[B]), de Jacobi et à l'idéal de Stickelberger. Par rapport à [AB], on ne donnera les démonstrations que là où le cas abélien se distingue vraiment du cas cyclotomique.

\subsection{Modules et nombres de Weil}

\newtheorem{Definition}{Définition}[subsection]
\newtheorem{Conjecture}{Conjecture}[subsection]

\begin{Definition}
On appelle module de Weil de $K_n$ et l'on note $\mathcal{W}_n$ l'ensemble des $f$ appartenant à  $\mathrm{Hom}_{\mathbb{Z}[G_n]}(I_n,K_n^{\times})$ tels que
\begin{enumerate}
	\item $\exists \beta (f) \in \mathbb{Z}[G_n]\  tel\  que\  \forall \alpha \in P_n ,\  f(\alpha \mathcal{O}_n) \equiv \alpha^{\beta(f)}\ mod \ \mu_{K_n}$
	\item $f(I_n)\subset \mu_{K_n}\mathcal{U}_n.$ 
\end{enumerate}
\end{Definition}

\begin{Definition}
On appelle module des nombres de Weil de $K_n$ et l'on note $\mathrm{Weil}_n$ 
$$\mathrm{Weil}_n := \{ f(\mathfrak{a}) ; f\in \mathcal{W}_n, \mathfrak{a} \in I_n \}.$$
\end{Definition} 

Notons que $\mathrm{Weil}_n$ est un sous-module de $\mu_{K_n}\mathcal{U}_n$.
Ce module ne servira qu'à partir de la deuxième partie. 

\newtheorem{Lemme}{Lemme}[subsection]

Le résultat suivant généralise [Iw1], pp 102-103, [AB] Proposition 1, et sa démonstration est exactement celle de [AB] Proposition 1: 

\newtheorem{Theoreme}{Théorème}[subsection]

\begin{Theoreme}
On a une suite exacte de $\mathbb{Z}[G_n]$-modules
$$0\longrightarrow \emph{Hom}_{\mathbb{Z}[G_n]}(I_n,\mu_{K_n}) \longrightarrow \mathcal{W}_n^{-} \stackrel{\beta}{\longrightarrow} \emph{Ann}_{\mathbb{Z}[G_n]} \mathcal{C}l (K_n)^{-}\longrightarrow B_{n} \longrightarrow 0$$
où $B_n$ est un 2-groupe fini abélien, et $\beta$ : $f\mapsto \beta (f)$. 
\end{Theoreme}

\subsection{Sommes de Gauss, modules de Jacobi}

Posons $E_n=\mathbb{Q}(\mu_{c_n})$ le corps conducteur de $K_n$. 
Soit $\mathcal{O}_{E_n}$ l'anneau des entiers de $E_n$ et $G_{E_n}$ le groupe de Galois de $E_n/\mathbb{Q}$. Soit $\mathfrak{L}$ un idéal premier de $E_n$ ne divisant pas $c_n$ et $(l)=\mathfrak{L}\cap\mathbb{Z}$. Notons $\mathcal{O}_{E_n} /\mathfrak{L}=\mathbb{F}_q$.
\\On définit le caractère $$\chi_{\mathfrak{L}} : \mathcal{O}_{E_n} /\mathfrak{L} \rightarrow \mu_{c_n}$$ pour $y\neq 0$ par la relation de congruence $$ \chi_{\mathfrak{L}} ( y \ \mathrm{mod}\  \mathfrak{L} ) \equiv y^{- \frac{q-1}{c_n}}\  \mathrm{mod} \ \mathfrak{L}$$
et l'on étend $\chi_{\mathfrak{L}}$ à $\mathbb{F}_q$ en posant $\chi_{\mathfrak{L}}(0)=0$.
\\Soit Tr la trace de $\mathbb{F}_q  / \ \mathbb{F}_l$. On a un caractère
$$\psi_q : \mathcal{O}_n /\mathfrak{L} \rightarrow \mu_l, \  \  \psi_q(x)= \zeta_l^{\mathrm{Tr}(x)}$$
où $\zeta_l$ est une racine primitive $l$-ième de l'unité.
\\On pose alors $$\tau_{E_n} (\mathfrak{L}) := - \displaystyle{\sum_{a \in \mathbb{F}_q}} \chi_{\mathfrak{L}}(a) \psi_q (a)$$
et l'on dit que $\tau_{E_n} (\mathfrak{L})$ est la \textit{somme de Gauss} associée à $\mathfrak{L}$.
\\Comme tous les idéaux au-dessus de $l$ ont le même corps résiduel associé à isomorphisme près, comme $G_{E_n}$ agit par permutation sur ceux-ci et comme $\chi_{\mathfrak{L}}^{\delta}=\chi_{\mathfrak{L}^{\delta}}$ pour tout $\delta \in G_{E_n}$, on a, pour tout $\mathfrak{L}$ idéal premier de $E_n$ ne divisant pas $c_0$,
$$\forall \delta \in G_{E_n}, \ \tau_{E_n}(\mathfrak{L})^{\delta}=\tau_{E_n}(\mathfrak{L^{\delta}}).$$
Donc $\tau_{E_n}$ s'étend multiplicativement en un morphisme de $\mathbb{Z}[G_{E_n}]$-modules
$$\tau_{E_n} : \ I_{E_n} \rightarrow \Omega(\mu_{c_n})^{*}$$
où $\Omega$ est le compositum de tous les $\mathbb{Q}(\mu_m)$, $m$ premier à $c_n$ et en notant $I_{E_n}$ l'ensemble de tous $\mathfrak{L}$ idéaux premiers de $E_n$ ne divisant pas $c_n$.
\\Par construction de $\Omega$, $G_n$ est canoniquement isomorphe à $Gal(\Omega K_n/\Omega)$ donc $\tau_{E_n}$ induit un morphisme de $\mathbb{Z}[G_n]$-modules
$$\tau_n : \ I_n \rightarrow \Omega(\mu_{c_n})^{*}.$$
  Soit $t$ l'ordre du groupe multiplicatif du corps résiduel d'un idéal premier de $K$ au-dessus de $p$. Notons que $t$ ne dépend pas du choix d'un tel idéal, puisque $G_0$ agit transitivement sur ceux-ci. Alors, pour tout $x\in K_n^{*}$ tel que $(x,p)=1$, on a $x^{t}\in \mathcal{U}_n$ via le plongement diagonal. Posons
$$\widetilde{\tau}_n := \tau_n^{t}.$$

\subsection{Sommes de Jacobi et théorème de Stickelberger}

\begin{Definition}
On appelle module de Jacobi de $K_n$ et l'on note $\mathcal{J}_n$ l'ensemble défini par
$$\mathcal{J}_n:=\mathbb{Z}[G_n] \widetilde{\tau}_n \cap \mathrm{Hom}_{\mathbb{Z}[G_n]}(I_n, K_n^{*}).$$
On appelle module des sommes de Jacobi de $K_n$ et l'on note $\mathrm{Jac}_n$ l'ensemble défini par
$$\mathrm{Jac}_n := \{ f(\mathfrak{a});\  f\in\mathcal{J}_n,\mathfrak{a}\in I_n \}.$$
\end{Definition}

Soit $a$ un entier premier à $c_n$. Soit $\sigma_a \in G_{E_n}$  l'image de $a$ mod $c_n \mathbb{Z}$ via l'isomorphisme
$$(\mathbb{Z}/c_n\mathbb{Z})^{\times} \simeq G_{E_n}.$$
L'élément de Stickelberger de $E_n$ est alors $$\Theta_{E_n} := \displaystyle{\sum_{a=1..c_n, (a,c_n)=1}}  \frac{a}{c_n}  \ \sigma_a^{-1}\in \mathbb{Q}[G_{E_n}].$$
Et celui de $K_n$ est la restriction de $\Theta_{E_n}$ à $K_n$ au sens des morphismes galoisiens, noté $\Theta_n$.

\begin{Definition}
L'idéal de Stickelberger de $K_n$ est
$$\mathcal{S}_n := \mathbb{Z}[G_n]\cap \Theta_n \mathbb{Z}[G_n].$$
\end{Definition}

Posons $\mathcal{S}'_n=\{g\in\mathbb{Z}[G_n]$ ; $g t\Theta_n\in\mathbb{Z}[G_n]\}$ où $t$ a été défini au paragraphe précédent. 

\begin{Theoreme}[Théorème de Stickelberger]
Soit $\mathfrak{p}$ un idéal premier de $I_n$. Supposons que $\beta\Theta_n\in\mathbb{Z}[G_n]$. Alors $\tau_n(\mathfrak{p})^{\beta} \in K_n^{\times}$, et
$$\tau_n(\mathfrak{p})^{\beta} \mathcal{O}_n = \mathfrak{p}^{\beta \Theta_n}.$$
Donc, si $\beta \in \mathcal{S}'_n$, $\widetilde{\tau}_n(\mathfrak{p})^{\beta} \in \mathcal{U}_n$, $\beta t \Theta_n \in \mathbb{Z}[G_n]$ et
$$\widetilde{\tau}_n(\mathfrak{p})^{\beta} \mathcal{O}_n = \mathfrak{p}^{\beta t \Theta_n}.$$
\end{Theoreme}

\textbf{Remarques.}
\\-- Notons que ce théorème est un résultat d'annulation du groupe de classes d'idéaux de $K_n$.
\\-- La première partie de la preuve de ce théorème (voir [Wa, p.98]) nous donne le résultat:
	\\Soit $l\equiv 1$ mod $c_n$ un nombre premier, $\mathfrak{L}$ au-dessus de $l$ dans $E_n$ et $\tilde{\mathfrak{L}}$ l'unique idéal au-dessus de $\mathfrak{L}$ dans $F_n(\mu_l)$. Alors $$\tau_n(\mathfrak{L})\mathbb{Z}[\zeta_{q_n},\zeta_l]=\tilde{\mathfrak{L}}^{(l-1)\Theta_{E_n}}.$$ Ce résultat implique sans difficulté le lemme suivant (voir [AB] Lemme2):

\begin{Lemme} 
$$\mathcal{S}'_n=\{ u \in \mathbb{Z}[G_n], \widetilde{\tau}_n^{u}\in\mathcal{J}_n \};$$
$$\mathcal{J}_n=\{ \widetilde{\tau}_n^{\delta}, \delta\in \mathcal{S}'_n \}.$$
\end{Lemme}

\newtheorem{Proposition}{Proposition}[subsection]

\begin{Proposition}
Le morphisme $\beta$ induit par restriction le diagramme commutatif:
\ \\
	\xymatrix{
& & &	&  \mathcal{W}_n \ar@{->}[r] & \mathrm{Ann}_{\mathbb{Z}[G_n]}(\mathcal{C}l(K_n))\\
& & & &  \mathcal{J}_n  \ar@{^(->>}[r] \ar@{^(->}[u]^{Id} & t\mathcal{S}_n \ar@{^(->}[u]^{Id}}
\end{Proposition}

\begin{proof}[\normalfont{\bfseries{Preuve.}}] -- Notons $N_n$ l'élément norme de $G_n$. Soit $f\in\mathcal{J}_n$. Par le lemme 2.3.1, il existe $\delta \in \mathcal{S}'_n$ tel que  $f=\widetilde{\tau}_n^{\delta}$. Soit $\alpha\mathcal{O}_{K_n} \in I_n$. Par le théorème de Stickelberger, on a $f(\alpha\mathcal{O}_{K_n})=\widetilde{\tau}_n^{\delta}(\alpha\mathcal{O}_{K_n})=\epsilon \alpha^{\delta t\Theta_n}$ avec $\epsilon \in \mathcal{O}_{K_n}^{\times}$. Mais [Wa, lemma 6.1], $\widetilde{\tau}_n^{\delta}(\alpha\mathcal{O}_{K_n})\widetilde{\tau}_n^{\delta\sigma_{-1}}(\alpha\mathcal{O}_{K_n})=\alpha^{t \delta N_n}.$ Donc, comme $G_n$ est commutatif, $\alpha^{t \delta N_n}=\epsilon \epsilon^{\sigma_{-1}} \alpha^{t \delta(\Theta_n+\sigma_{-1}\Theta_n)}.$ Notons $\{x\}$ la partie fractionnaire du nombre réel $x$. Comme $\{x\}+\{-x\}=1$ pour tout $x\in\mathbb{Q}\backslash\mathbb{Z}$, $N_n=\Theta_n+\sigma_{-1}\Theta_n$. Donc $\epsilon \ \epsilon^{\sigma_{-1}}=1$. D'où $\epsilon\in\mu_{q_n}$. Finalement,	
\begin{equation}
f(\alpha\mathcal{O}_{K_n})\equiv \alpha^{\delta t \Theta_n} \ mod \ \mu_{q_n}.
\end{equation}
	Mais le point 2 de la définition 2.1.1 est vérifié par définition de $\widetilde{\tau}_n$, donc $f\in \mathcal{W}_n$.
	\\-- On a un diagramme commutatif:
	\\ \ \\
	\xymatrix{
&  \mathrm{Hom}_{\mathbb{Z}[G_n]}(I_n,\mu_{q_n})  \ar@{^(->}[r] &  \mathcal{W}_n \ar@{->}[r] & \mathrm{Ann}_{\mathbb{Z}[G_n]}(\mathcal{C}l(K_n))\\
& \mathrm{Ker} \beta' \ar@{^(->}[r]
\ar@{^(->}[u]& \mathcal{J}_n  \ar@{-->>}[r]^{\beta'} \ar@{^(->}[u] & t\mathcal{S}_n \ar@{^(->}[u].}
\\ \ \\ \ \\
La ligne du dessus est donnée par le théorème 2.1.1. Les inclusions de $\mathcal{J}_n$ et $\mathcal{S}_n$ proviennent de 1. et du théorème de Stickelberger respectivement. De plus par (1), pour tout $f\in\mathcal{J}_n$, $f=\widetilde{\tau}_n^{\delta}$, $\delta\in\mathcal{S}'_n$, $\beta(f)=\beta(\widetilde{\tau}_n^{\delta})=\delta t \Theta_n.$ Donc $\beta$ induit par restriction à $\mathcal{J}_n$ une application $\beta'$, surjective par le lemme 2.3.1. et la définition de $\mathcal{S}'_n$. Soit $\delta\in\mathcal{S}'_n$ tel que $\delta t \Theta_n=0$. On a 
$\mathrm{Ker}(\beta')\subset \mathrm{Hom}_{\mathbb{Z}[G_n]}(I_n,\mu_{q_n}) \subset \mathcal{W}_n^{-}$, donc $\widetilde{\tau}_n^{\sigma_{-1} \delta}=\widetilde{\tau}_n^{-\delta}$. En outre, l'élément de Stickelberger n'est pas un diviseur de zéro dans $\mathbb{Z}[G_n]^{-}$, donc $\delta t \in\mathbb{Z}[G_n]^{+}$. Donc $\widetilde{\tau}_n^{2\delta}=1$. Mais $\widetilde{\tau}_n^{\delta}\in\mathcal{U}_n$, c'est-à-dire $\widetilde{\tau}_n^{\delta}=1$.

\end{proof}

Ce diagramme commutatif peut suggérer une extension de la définition de $\mathcal{J}_n$ au cas relativement abélien, c'est-à-dire où le corps de base n'est pas $\mathbb{Q}$ (celles de $\mathcal{W}_n$ et $\mathrm{Weil}_n$ s'étendant sans difficulté). Plus précisément, prenons (juste pour cette fin de paragraphe) les notations suivantes: 
\\$F$ est un corps de nombres totalement réel. Soit $K/F$ une extension abélienne de $F$. Supposons que $\mu_p\subset K$.
\\Notons, pour $n\geq1$ $K_n=K(\mu_{p^{n+1}})$, $G_n=\mathrm{Gal}[K_n/F]$.

\textsl{L'élément de Stickelberger} de $K_n$ est défini par
$$\Theta_n:=\displaystyle{\sum_{\sigma\in G_n}}\zeta_{F,R_n}(\sigma,0)\sigma^{-1},$$
où $\zeta_{F,R_n}(\sigma,s)$ est la fonction zeta partielle, définie pour $Re(s)>1$ par
$$\zeta_{F,R_n}(\sigma,s)=\displaystyle{\sum_{(\mathfrak{a},R_n)=1, (\mathfrak{a},K_n/F)=\sigma}}N\mathfrak{a}^{-s}.$$
et admettant un prolongement méromorphe à $\mathbb{C}$.
\\On sait que l'ensemble
$$\mathcal{S}_n:=Ann_{\mathbb{Z}[G_n]}(\mu_{K_n})\Theta_n$$
est un idéal de l'anneau $\mathbb{Z}[G_n]$. La conjecture de Brumer s'énonce alors

\begin{Conjecture}
$\mathcal{S}_n^{-}\subset(\mathrm{Ann}_{\mathbb{Z}[G_n]}\mathcal{C}l_{K_n})^{-}$.
\end{Conjecture}

Par analogie avec le cas cyclotomique, on appelle \textsl{module de Jacobi} et l'on note $\mathcal{J}_n$ l'image réciproque de $t\mathcal{S}_n$ par $\beta$. On appelle \textsl{module des sommes de Jacobi} et l'on note $\mathrm{Jac}_n$ le module défini par 
$$\mathrm{Jac}_n:=\{f(\mathfrak{p}); \ f\in \mathcal{J}_n, \ \mathfrak{p}\in I_n\}.$$
Avec ces définitions et la conjecture de Brumer, on pourrait généraliser une partie des résultats de cet article.

\subsection{$\psi$-parties, quelques propriétés fonctorielles}

Pour tout $\mathbb{Z}$-module $M$, notons $\overline{M}$ le pro-$p$-complété de $M$. On a $\overline{M}=M\otimes\mathbb{Z}_p$ si $M$ est de type fini. Comme $\mathbb{Z}_p$ est plat et $p\neq2$, la suite exacte du théorème 2.1.1. donne
$$0\longrightarrow \mathrm{Hom}_{\mathbb{Z}[G_n]}(I_n,\mu_{K_n})\otimes_{\mathbb{Z}}\mathbb{Z}_p \longrightarrow \overline{\mathcal{W}}_n^{-} \stackrel{\overline{\beta}}{\longrightarrow} (\mathrm{Ann}_{\mathbb{Z}_p[G_n]} A_n)^{-} \longrightarrow 0,$$
où $A_n$ est le $p$-Sylow de $\mathcal{C}l_{K_n}$ et $\overline{\beta}$ est naturellement induit par $\beta$. 
\\ \ \\
\textsl{A partir d'ici et sauf mention expresse du contraire, tous les calculs caractère par caractère se feront sous l'hypothèse de semi-simplicité, i.e. $p$ ne divise pas} [K:$\mathbb{Q}$].
\\Soit $\psi$ un $\overline{\mathbb{Q}}_p$ caractère impair irréductible de $\Delta$, différent du caractère de Teichmüller $\omega$. Notons $\psi^{*}:=\omega.\psi^{-1}$ son caractère "miroir", $\mathbb{Z}_p[\psi]$ le sous anneau de $\overline{\mathbb{Q}}_p$ engendré par l'image de $\psi$. L'idempotent  $e_{\psi}$ associé à $\psi$ est défini par
$$e_{\psi}=\frac{1}{|\Delta|}\displaystyle{\sum_{\delta\in\Delta}}\psi(\delta)\delta^{-1}\in e_{\psi}\mathbb{Z}_p[\Delta]\simeq\mathbb{Z}_p[\psi].$$
Soit $M$ un $\mathbb{Z}_p[\Delta]$-module. Dans toute la suite, on notera $M^{\psi}=e_{\psi}M$.
Soit $\Lambda_{\psi}=\displaystyle{\lim_{\longleftarrow}}\ e_{\psi}\mathbb{Z}_p[G_n]=\displaystyle{\lim_{\longleftarrow}}\    e_{\psi}\mathbb{Z}_p[\Delta][\Gamma_n]\simeq\mathbb{Z}_p[\psi][[\Gamma]]$. 
\\Soit $\gamma_0$ un générateur topologique de $\Gamma$. On a l'isomorphisme bien connu $\Lambda_{\psi}\simeq\mathbb{Z}_p[\psi][[T]]$, en envoyant $\gamma_0-1$ sur $T$.
\\Comme $\psi\neq\omega$ et $t$ (défini au §2.2) est premier à $p$, le diagramme du paragraphe précédent induit le diagramme commutatif de $\mathbb{Z}_p[\psi][\Gamma_n]$-modules:
\\
\xymatrix{ & & & & \overline{\mathcal{W}}_{n}^{\psi} \ar@{^(->>}[r] & \mathrm{Ann}_{\mathbb{Z}_p[\psi][\Gamma_n]} A_n \\ 
& & & & \overline{\mathcal{J}}_{n}^{\psi} \ar@{^(->>}[r] \ar@{^(->}[u] & \overline{\mathcal{S}}_{n}^{\psi}. \ar@{^(->}[u]
}
\ \\
Nous allons maintenant montrer que les objets de ce diagramme "montent".
Comme $\psi\neq\omega$, $\overline{\mathcal{S}}_{n}^{\psi}$ est principal, engendré par $\Theta_{n}^{\psi}$ la $\psi$-partie de l'élément de Stickelberger de $K_n$, qui est alors un élément de $\mathbb{Z}_p[\psi][G_n]$.
\\De plus, $(\Theta_{n}^{\psi})_{n\geq0}$ forme un système projectif pour la restriction des automorphismes (voir e.g. [Wa, Proposition 7.6]) . Notons
$$\Theta_{\psi}:=\displaystyle{\lim_{\longleftarrow}}\ \Theta_{n}^{\psi}$$ 
et $X_{\infty}=\displaystyle{\lim_{\longleftarrow}}\  A_n$, la limite étant prise pour les applications de norme. Il est clair que
$$\displaystyle{\lim_{\longleftarrow}}\ \mathrm{Ann}_{\mathbb{Z}_p[\psi][\Gamma_n]}A_{n}^{\psi}= \mathrm{Ann}_{\Lambda_{\psi}} X_{\infty}^{\psi}.$$
Les objets de gauche du diagramme montent alors via l'isomorphisme $\overline{\beta}$. Notons leur limites respectives $$\overline{\mathcal{J}}_{\infty}^{\psi}, \ \overline{\mathcal{W}}_{\infty}^{\psi}.$$
Comme $(\mathcal{U}_n)_{n\geq 0}$ et $(I_n)_{n\geq 0}$ sont évidemment des systèmes projectifs pour la norme, $(\overline{\mathrm{Weil}}_{n}^{\psi})_{n\geq 0}$ et $(\overline{\mathrm{Jac}}_{n}^{\psi})_{n\geq 0}$ en sont aussi.  
Notons leur limites respectives
$$\mathcal{U}_{\infty}
,\ \overline{\mathrm{Weil}}_{\infty}^{\psi},\ \overline{\mathrm{Jac}}_{\infty}^{\psi}.$$

\section{Montée et annulateurs}
En profitant des propriétés fonctorielles précédentes, on va maintenant passer à la limite projective le long de la tour cyclotomique pour paramétrer, au niveau infini, les annulateurs du $\Lambda_{\psi}$-module $X_{\infty}^{\psi}$.

\subsection{Idéaux annulateurs}

Comme depuis le début, $\psi$ est impair, $\psi\neq\omega$.
Le lemme suivant ([H, p.31, remarque 2]-[Wa, proposition 13.28]) permet de définir la "série minimale" d'un certain type de module de torsion:

\begin{Lemme}
Si $M$ est de $\Lambda_{\psi}$-torsion, sans sous-module fini non nul, alors $M$ est de dimension projective inférieure ou égale à un sur $\Lambda_{\psi}$ (Auslander-Buchsbaum), et par suite $M$ admet une matrice de présentation $P$ qui est carrée, disons de dimension $n$. Soit $\Delta_0=\mathrm{det}(P)$ et $\Delta_1$ le plus grand diviseur commun des $(n-1)$-mineurs de $P$. L'idéal de Fitting et l'annulateur de $M$ sont principaux sur $\Lambda_{\psi}$, engendrés resp. par $\Delta_0$ et $\Delta_0/\Delta_1$.
\end{Lemme}

Sous les hypothèses du lemme 3.1.1., on appelle \textsl{série minimale} de $M$ un générateur de l'annulateur de $M$. C'est un annulateur dont le degré de Weierstrass est minimal. Bien entendu, la série minimale de $M$ divise sa série caractéristique.
\\Si $E$ est le module élémentaire associé à $M$, $M$ et $E$ ont même annulateur. La matrice de présentation de $E$ étant diagonale, un générateur de l'annulateur est le plus petit commun multiple des termes diagonaux (i.e. des termes qui apparaissent dans le théorème de structure) alors que la série caractéristique en est le produit. 
\\ \ \\
\textbf{Fait:} $X_{\infty}^{\psi}$ vérifie ces hypothèses pour $\psi$ impair (voir e.g [W, §13]).

\subsection{Quelques représentations galoisiennes}

Dans cette section, et sauf mention expresse du contraire, $\psi$ est un caractère impair, \textsl{distinct du caractère de Teichmüller $\omega$}.
\\On se propose de construire certaines représentations galoisiennes de $\mathfrak{X}_n^{\psi}$ (resp. $\mathfrak{X}_{\infty}^{\psi}$) dans $\mathbb{Z}_p[\psi][\Gamma_n]$ (resp. $\Lambda_{\psi}$) qui joueront un rôle important dans l'étude des annulateurs galoisiens du groupe de classes.
\\Notre premier résultat va décrire les modules $\overline{\mathrm{Jac}}_n^{\psi}$ et $\overline{\mathrm{Jac}}_{\infty}^{\psi}$ comme des modules de sommes de Gauss au sens d'Ichimura-Hachimori-Beliaeva ([HI], §3.2; [B], §6).
\ \\
On va utiliser les notations suivantes:
\\$L_n$ la pro-$p$-extension abélienne non ramifiée maximale de $K_n$.
\\$L'_n$ la pro-$p$-extension abélienne non ramifiée p décomposée maximale de $K_n$.
\\$M_n$ la pro-$p$-extension abélienne non ramifiée hors de $p$ maximale de $K_n$.
\\$X_n=\mathrm{Gal}(L_n/K_n)$, $X'_n=\mathrm{Gal}(L'_n/K_n)$ et enfin $\mathfrak{X}_n=\mathrm{Gal}(M_n/K_n)$.
Nous noterons par un indice infini les limites projectives des groupes de Galois pour les applications de norme.

\begin{Lemme}
Tout $z_n\in\overline{\mathcal{W}}_n^{\psi}$ (resp. $z_{\infty}\in\overline{\mathcal{W}}_{\infty}^{\psi}$) induit un morphisme de $\mathbb{Z}_p[\psi][\Gamma_n]$-modules $z_n:\ \mathfrak{X}_n^{\psi}\rightarrow \mathcal{U}_n^{\psi}$ (resp. un morphisme de $\Lambda_{\psi}$-modules  $z_{\infty}:\ \mathfrak{X}_{\infty}^{\psi}\rightarrow \mathcal{U}_{\infty}^{\psi}$). Si l'on identifie par le corps de classes $\mathcal{U}_n^{\psi}$ (resp. $\mathcal{U}_{\infty}^{\psi}$) à un sous-module de $\mathfrak{X}_n^{\psi}$ (resp. de $\mathfrak{X}_{\infty}^{\psi}$), $z_n$ (resp. $z_{\infty}$) s'identifie à la multiplication par $\beta(z_n)$ (resp. $\beta(z_{\infty})$) dans $\mathfrak{X}_n^{\psi}$ (resp. $\mathfrak{X}_{\infty}^{\psi}$).  
\end{Lemme}

\begin{proof}[\normalfont{\bfseries{Preuve.}}]
Comme $\mathcal{U}_n$ est un $\mathbb{Z}_p$-module, le morphisme $z_n$ induit par la théorie du corps de classes un morphisme, encore noté $z_n$, de $\mathbb{Z}_p[\psi][\Gamma_n]$-modules
$$z_n:\ \mathfrak{X}_{n}^{\psi} \rightarrow \mathcal{U}_{n}^{\psi}.$$ 
Si un système $(z_n)_{n\geq 0}$ est cohérent pour la norme, il définit un homomorphisme de $\Lambda_{\psi}$-modules:

$$z_{\infty}:\ \mathfrak{X}_{\infty}^{\psi}\rightarrow\mathcal{U}_{\infty}^{\psi}.$$
\ \\L'identification de $z_n$ (resp. $z_{\infty}$) résulte du Théorème 2.1.1. et des suites exactes du corps de classes ($\psi\neq \omega$, $\psi$ impair):
$$0\rightarrow\mathcal{U}_n^{\psi}\rightarrow\mathfrak{X}_n^{\psi}\rightarrow A_n^{\psi}\rightarrow 0.$$

\end{proof}

A partir de maintenant, on notera $\mathcal{F}_{\psi}$ "la" série caractéristique de $X_{\infty}^{\psi}$, $\mathcal{M}_{\psi}$ "sa" série minimale au sens du lemme 3.1.1. D'après la Conjecture Principale (ou théorème de Mazur-Wiles),  les séries $\Theta_{\psi}$ et $\mathcal{F}_{\psi}$ diffèrent par multiplication par un élément inversible de $\Lambda_{\psi}$. En particulier, $\Theta_{\psi}\Lambda_{\psi}=\mathcal{F}_{\psi}\Lambda_{\psi}$.

\begin{Lemme}
$\overline{\mathrm{Jac}}_n^{\psi}=\Theta_n^{\psi}\mathfrak{X}_n^{\psi}$ et $\overline{\mathrm{Jac}}_{\infty}^{\psi}=\Theta_{\psi}\mathfrak{X}_{\infty}^{\psi}=\mathcal{F}_{\psi}\mathfrak{X}_{\infty}^{\psi}$ 
\end{Lemme}
C'est la caractérisation des \textsl{modules de sommes de Gauss} $p$-adiques au niveau infini par Hachimori-Ichimura-Beliaeva ([HI], [B]).
\\Remarquons au passage que comme $\overline{\mathcal{W}}_{\infty}^{\psi}\simeq \mathcal{M}_{\psi}\Lambda_{\psi}$, $\overline{\mathrm{Weil}}_{\infty}^{\psi}=\mathcal{M}_{\psi}\mathfrak{X}_{\infty}^{\psi}$.

\begin{proof}[\normalfont{\bfseries{Preuve du lemme 3.2.2.}}]
D'après le lemme 2.3.1., tout élément de $\mathcal{J}_n$ est de la forme $\widetilde{\tau}_n^{\delta}$, donc $\overline{\mathrm{Jac}}_n^{\psi}\subset \mathrm{Im} \widetilde{\tau}_n$, et d'après le lemme 2.3.1., $\mathrm{Im} \widetilde{\tau}_n=\beta(\widetilde{\tau}_n)\mathfrak{X}_n^{\psi}$. Or $\beta(\widetilde{\tau}_n)=t\Theta_n$ et $t$ est premier à $p$.
\end{proof}

Dans la suite, on verra que la $\Lambda_{\psi}$-liberté éventuelle de $\mathcal{U}_{\infty}^{\psi}$ joue un rôle central. Si $\psi^*(p)\neq 1$, c'est un module libre [G, proposition 1], mais dans le cas $\psi^*(p)=1$, on va devoir modifier les modules $\mathcal{U}_{\infty}^{\psi}$, $\mathfrak{X}_{\infty}^{\psi}$ et $\mathrm{Jac}_{\infty}^{\psi}$ comme suit:
\\Notons $\mathfrak{T}_n:=\mathrm{Tor}_{\mathbb{Z}_p}\mathcal{U}_n$ et $\mathfrak{T}=\displaystyle{\lim_{\longleftarrow}}\ \mathfrak{T}_n$, et remarquons que $\mathfrak{T}_n^{\psi}=\{1\}$ si $\psi^{*}(p)\neq1$. 
Notons $\widetilde{\mathcal{U}}_n^{\psi}:=\mathcal{U}_n^{\psi}/\mathfrak{T}_n^{\psi}$ (resp. $\widetilde{\mathcal{U}}_{\infty}^{\psi}:=\mathcal{U}_{\infty}^{\psi}/\mathfrak{T}^{\psi}$) et $\widetilde{\mathfrak{X}}_n^{\psi}:=\mathfrak{X}_n^{\psi}/\mathfrak{T}_n^{\psi}$ (resp. $\widetilde{\mathfrak{X}}_{\infty}^{\psi}:=\mathfrak{X}_{\infty}^{\psi}/\mathfrak{T}^{\psi}$) en voyant $\mathcal{U}_n^{\psi}$ comme un sous-module de $\mathfrak{X}_n^{\psi}$ (resp. $\mathcal{U}_{\infty}^{\psi}$ comme un sous-module de $\mathfrak{X}_{\infty}^{\psi}$) via le corps de classes. Le module $\widetilde{\mathfrak{X}}_{\infty}^{\psi}$ est la $\psi$-partie du \textbf{module de Bertrandias-Payan} dont les propriétés sont rappelées dans l'appendice.
\\Les applications $z_n$ (resp. $z_{\infty}$) se factorisent alors en des morphismes de $\mathbb{Z}_p[\psi][\Gamma_n]$-modules (resp. de $\Lambda_{\psi}$-modules)
$$z_n:\ \widetilde{\mathfrak{X}}_n^{\psi}\rightarrow \widetilde{\mathcal{U}}_n^{\psi},\  z_{\infty}:\ \widetilde{\mathfrak{X}}_{\infty}^{\psi}\rightarrow \widetilde{\mathcal{U}}_{\infty}^{\psi}.$$
\\On sait (voir e.g. [G], propositions 1 et 2) que $\widetilde{\mathcal{U}}_{\infty}^{\psi}\simeq\Lambda_{\psi}$ et que si $\psi(p)\neq1$, $\widetilde{\mathcal{U}}_n^{\psi}\simeq\Lambda_{\psi}/\omega_n\Lambda_{\psi}\simeq\mathbb{Z}_p[\psi][\Gamma_n]$ ($\omega_n=((1+T)^{p^n}-1)$), donc les morphismes $z_{\infty}$ et $z_n$ si $\psi(p)\neq1$ peuvent être considérés, ce qu'on fera désormais, comme des représentations galoisiennes
$$z_n:\ \widetilde{\mathfrak{X}}_n^{\psi}\rightarrow\mathbb{Z}_p[\psi][\Gamma_n], \ z_{\infty}:\ \widetilde{\mathfrak{X}}_{\infty}^{\psi}\rightarrow\Lambda_{\psi}.$$

\subsection{Une suite exacte fondamentale}

On va maintenant, à partir de $X_{\infty}^{\psi}$, fabriquer une suite exacte qui donnera des renseignements sur le "miroir" de $X_{\infty}^{'\psi}$ (défini au § 3.2).

Soit $Z\in\Lambda_{\psi}$ tel que $\mathcal{M}_{\psi}$ divise $Z$ et $Z$ divise $\mathcal{F}_{\psi}$ dans $\Lambda_{\psi}$.
Soit $(Z_n)_{n\in\mathbb{N}}$ l'image de $Z$ par la surjection canonique $\Lambda_{\psi}\rightarrow\Lambda_{\psi}/\omega_n$. On a 
$$Z_n\equiv Z\ \mathrm{mod} \ \omega_n,\ Z=\displaystyle{\lim_{\longleftarrow}}\ Z_n.$$

\begin{Lemme}
Pour n>>0, $Z_n\in \mathrm{Ann}_{\mathbb{Z}_p[\psi] [\Gamma_n]} A_n^{\psi}$.
\end{Lemme}

\begin{proof}[\normalfont{\bfseries{Preuve.}}]
Par définition, $Z$ annule $X_{\infty}^{\psi}$, donc $Z_n$ annule $X_{\infty}^{\psi}/\omega_n X_{\infty}^{\psi}$. Or ce dernier quotient se surjecte naturellement sur $A_n^{\psi}$ pour $n$>>$0$.
\end{proof}

Construisons $z_n$, $z_{\infty}$ par la méthode du §3.2 à partir des images des $Z_n$ par les isomorphismes $\overline{\mathcal{W}}_n^{\psi}\simeq \mathrm{Ann}_{\mathbb{Z}_p[\psi] [\Gamma_n]} A_n^{\psi}$.
\begin{equation}
z_n:\ \widetilde{\mathfrak{X}}_n^{\psi}\rightarrow\mathbb{Z}_p[\psi][\Gamma_n],\ z_{\infty}:\ \widetilde{\mathfrak{X}}_{\infty}^{\psi}\rightarrow \Lambda_{\psi}.
\end{equation}
Dans la suite, les applications correspondant à $\mathcal{M}_{\psi}$ et $\mathcal{F}_{\psi}$ joueront un rôle particulier, on les notera $m_n$, $m_{\infty}$ resp. $f_n$, $f_{\infty}$.
Notons que si cette construction est licite pour $n$>>0 par le lemme 3.3.1., elle l'est en général dès que $Z_n\in\mathrm{Ann}_{\mathbb{Z}_p[\psi] [\Gamma_n]} A_n^{\psi}$. En particulier, $f_n$ est définie pout tout $n$.

\begin{Theoreme} Soit $\psi$ impair, $\psi\neq\omega$.
\\L'application $z_{\infty}$ définie en (2) donne une suite exacte de $\Lambda_{\psi}$-modules:
$$0\rightarrow \mathrm{tor}_{\Lambda_{\psi}}\widetilde{\mathfrak{X}}_{\infty}^{\psi} \rightarrow X_{\infty}^{\psi} \rightarrow \Lambda_{\psi}/(Z)\rightarrow \Lambda_{\psi}/z_{\infty}(\widetilde{\mathfrak{X}}_{\infty}^{\psi})\rightarrow 0.$$
Au niveau fini, elle donne, si $\psi(p)\neq 1$, une suite exacte de $\mathbb{Z}_p[\psi][\Gamma_n]$-modules:
$$0\rightarrow \mathrm{tor}_{\mathbb{Z}_p}\widetilde{\mathfrak{X}}_n^{\psi} \rightarrow A_n^{\psi} \rightarrow\mathbb{Z}_p[\psi][\Gamma_n]/(Z_n)\rightarrow\mathbb{Z}_p[\psi][\Gamma_n]/z_n(\widetilde{\mathfrak{X}}_n^{\psi})\rightarrow 0.$$
\end{Theoreme}

\textbf{Commentaires.} Les suites exactes du théorème 3.3.1., au niveau fini comme au niveau infini, sont implicitement contenues dans les suites exactes du serpent (à six termes) qui apparaissent dans la démonstration des théorèmes 1 et 2 de [B]. Au niveau infini (théorème 1), notre module $X_{\infty}^{\psi}$ est remplacé par un sous-module (noté $D$) car [B] ne fait pas d'hypothèse de semi-simplicité; au niveau fini (théorème 2), la démonstration de la p.150 de [B] introduit l'hypothèse supplémentaire $\psi^{*}(p)\neq 1$.

\begin{proof}[\normalfont{\bfseries{Preuve du théorème 3.3.1.}}]
La preuve de ce théorème au niveau infini est essentiellement contenue dans [B, p .140]; prouvons-le au niveau fini. Il s'agit d'appliquer $z_n$ puis le lemme du serpent à la suite exacte du corps de classes relative à l'inertie:
$$\xymatrix{ & 0 \ar[r] & 0 \ar[r] \ar[d]& \mathrm{Ker}(z_n) \ar[r] \ar[d] & A_n^{\psi} \ar[d] \ar@{.>}`[rrdd]`[ldd]`[llllddd]`[rdddd][lldddd] & & \\
& 0 \ar[r] & \widetilde{\mathcal{U}}_n^{\psi} \ar[r] \ar[dd]^-{\times Z_n} & \widetilde{\mathfrak{X}}_n^{\psi} \ar[r] \ar[dd]^-{z_n}& A_n^{\psi} \ar[r] \ar[dd]^-{\alpha_n} & 0 &\\
& & & & & &\\
& 0 \ar[r] & \widetilde{\mathcal{U}}_n^{\psi} \ar[r] \ar[d] & \widetilde{\mathfrak{X}}_n^{\psi} \ar[r] \ar[d] & A_n^{\psi} \ar[r] \ar[d] & 0 &\\
& &  \widetilde{\mathcal{U}}_n^{\psi}/Z_n \widetilde{\mathcal{U}}_n^{\psi} \ar[r] & \widetilde{\mathfrak{X}}_n^{\psi}/z_n(\widetilde{\mathfrak{X}}_n^{\psi}) \ar[r] & A_n^{\psi} \ar[r] & 0 &}$$
La colonne de droite est claire: $\alpha_n$ est l'application nulle. Comme $\psi(p)\neq 1$, $\widetilde{\mathcal{U}}_n^{\psi}\simeq \widetilde{\mathcal{U}}_{\infty}^{\psi}/\omega_n \widetilde{\mathcal{U}}_{\infty}^{\psi} \simeq \Lambda_{\psi}/\omega_n \Lambda_{\psi} \simeq \mathbb{Z}_p[\psi][\Gamma_n]$ et par Gross, $Z_n$, qui est le co-descendu de $Z$, est premier à $\omega_n$. La flèche de gauche est donc injective car $\mathbb{Z}_p[\psi][\Gamma_n]$ n'a pas de sous-module fini non nul.
\\Reste à déterminer $\mathrm{Ker}(z_n)$. Par ce qu'on vient de montrer, ce module s'injecte dans $A_n^{\psi}$, qui est de $\mathbb{Z}_p$-torsion, donc $\mathrm{Ker}(z_n)\subset\mathrm{tor}_{\mathbb{Z}_p}(\widetilde{\mathfrak{X}}_n^{\psi})$. Pour montrer l'égalité, remarquons que la suite exacte d'inertie du corps de classes nous donne, comme $\widetilde{\mathcal{U}}_n^{\psi}$ est sans $\mathbb{Z}_p$-torsion, que $\mathrm{tor}_{\mathbb{Z}_p}(\widetilde{\mathfrak{X}}_n^{\psi})$ s'injecte dans $\mathrm{Ker}(z_n)$. Donc, comme $Z_n$ tue $A_n^{\psi}$, $z_n$ tue $\mathrm{tor}_{\mathbb{Z}_p}(\widetilde{\mathfrak{X}}_n^{\psi})$. On a ainsi la suite exacte:
$$0\rightarrow \mathrm{tor}_{\mathbb{Z}_p}(\widetilde{\mathfrak{X}}_n^{\psi}) \rightarrow A_n^{\psi} \rightarrow \widetilde{\mathcal{U}}_n^{\psi}/Z_n \widetilde{\mathcal{U}}_n^{\psi} \rightarrow \widetilde{\mathfrak{X}}_n^{\psi}/z_n(\widetilde{\mathfrak{X}}_n^{\psi}) \rightarrow A_n^{\psi} \rightarrow 0.$$
Le théorème suit en fixant un isomorphisme $\widetilde{\mathcal{U}}_n^{\psi}\simeq \mathbb{Z}_p[\Gamma_n]$, en notant que $z_n(\widetilde{\mathfrak{X}}_n^{\psi})\hookrightarrow \widetilde{\mathcal{U}}_n^{\psi}$ et que  $A_n^{\psi}\simeq \widetilde{\mathfrak{X}}_n^{\psi}/\widetilde{\mathcal{U}}_n^{\psi}$.
\end{proof}

\textbf{Remarque:}
\\Si $\psi(p)= 1$, non seulement $\widetilde{\mathcal{U}}_{\infty}^{\psi}/\omega_n$ n'est plus isomorphe à $\widetilde{\mathcal{U}}_n^{\psi}$ (voir e.g. [G]), mais $(Z,\omega_n)\neq 1$ et $\mathrm{Ker}(\times Z_n)$ n'est plus fini. On reviendra là-dessus dans la démonstration du thm. 3.4.2. ci-après.

Notons $\widetilde{\mathrm{Jac}}_{\infty}^{\psi}=\Theta_{\psi}\widetilde{\mathfrak{X}}_{\infty}^{\psi}=\mathcal{F}_{\psi}\widetilde{\mathfrak{X}}_{\infty}^{\psi}$ (resp. $\widetilde{\mathrm{Jac}}_{n}^{\psi}=\Theta_n^{\psi}\widetilde{\mathfrak{X}}_{n}^{\psi}$) et $\widetilde{\mathrm{Weil}}_{\infty}^{\psi}=\mathcal{M}_{\psi}\widetilde{\mathfrak{X}}_{\infty}^{\psi}$ (et remarquons que ces notations sont cohérentes par le lemme 3.2.2.).

De la démonstration du théorème on extrait le 

\newtheorem{Corollaire}{Corollaire}

\begin{Corollaire}
Avec les hypothèses et notations du théorème 3.3.1., on a des suites exactes compatibles par co-descente:
$$0 \longrightarrow \mathrm{tor}_{\Lambda_{\psi}} \widetilde{\mathfrak{X}}_{\infty}^{\psi} \longrightarrow \widetilde{\mathfrak{X}}_{\infty}^{\psi} \stackrel{\Theta_{\psi}}{\longrightarrow} \widetilde{\mathrm{Jac}}_{\infty}^{\psi}\simeq \mathrm{fr}_{\Lambda_{\psi}} \mathfrak{X}_{\infty}^{\psi} \longrightarrow 0$$
et 
$$0 \longrightarrow \mathrm{tor}_{\mathbb{Z}_p} \widetilde{\mathfrak{X}}_n^{\psi} \longrightarrow \widetilde{\mathfrak{X}}_n^{\psi} \stackrel{\Theta_n^{\psi}}{\longrightarrow} \widetilde{\mathrm{Jac}}_n^{\psi} \simeq \mathrm{fr}_{\mathbb{Z}_p} \mathfrak{X}_n^{\psi}\longrightarrow 0$$
où $\mathrm{fr}(.)$ désigne le quotient sans torsion.
\end{Corollaire}
  
Le théorème 3.3.1. donne aussi immédiatement "l'écart" entre l'idéal de Fitting et l'annulateur de $X_{\infty}^{\psi}$:

\begin{Corollaire}
Soient $Z,\ Z'\in\Lambda_{\psi}$ tels que $\mathcal{M}_{\psi}|Z|Z'|\mathcal{F}_{\psi}$. Alors
$$\frac{z_{\infty}(\widetilde{\mathfrak{X}}_{\infty}^{\psi})}{z'_{\infty}(\widetilde{\mathfrak{X}}_{\infty}^{\psi})}\simeq\frac{Z\Lambda_{\psi}}{Z'\Lambda_{\psi}}.$$
En particulier, $\mathcal{F}_{\psi}/\mathcal{M}_{\psi}=\mathrm{char}_{\Lambda_{\psi}}\widetilde{\mathrm{Weil}}_{\infty}^{\psi}/\widetilde{\mathrm{Jac}}_{\infty}^{\psi}$.
\end{Corollaire}

Comme dans [AB], on peut en tirer diverses formulations de versions plus ou moins affaiblies de la conjecture de Greenberg.

\subsection{Idéaux de Fitting}

Dans cette sous-section, on va déterminer les idéaux de Fitting initiaux des modules $X_{\infty}^{\psi^*}$ et $A_n^{\psi^*}$, pour \textbf{$\psi$ impair}, \textbf{distinct de $\omega$} , $\psi^{*}=\omega \psi^{-1}$ (on sait, par le théorème de Stickelberger que $A_n^{\omega}=\{0\}$). Les résultats seront énoncés en termes de modules de sommes de Gauss, mais notons qu'au niveau infini, ces modules interviennent uniquement sous la forme donnée dans le lemme 3.2.2. Pour ne pas se perdre dans les détails techniques, on sera amené parfois à supposer que $\psi(p)\neq 1$. En effet, le cas où $\psi(p)=1$ (i.e. $p$ est totalement décomposé dans la sous-extension fixée par $\mathrm{Ker}\psi$) est le cas notoirement le plus compliqué de la descente d'Iwasawa, celui où la nullité de certains homomorphismes naturels de descente empêche d'obtenir des renseignements au niveau fini (voir e.g. les commentaires dans l'introduction de [BelN]). Il faut alors modifier certains modules de $p$-unités circulaires, comme c'est fait par exemple dans [BelN], mais les complications techniques sont assez considérables pour qu'on ait préféré ici écarter le cas $\psi(p)=1$ - quitte à y revenir ultérieurement en cas de besoin.

\begin{Lemme}
Choisissons des isomorphismes cohérents $\widetilde{\mathcal{U}}_{\infty}^{\psi}\simeq\Lambda_{\psi}$ et $\widetilde{\mathcal{U}}_n^{\psi}\simeq\mathbb{Z}_p[\psi][\Gamma_n]$ (si $\psi(p)\neq1$), et notons $\widetilde{\mathcal{G}}_n^{\psi}$ et $\widetilde{\mathcal{G}}_{\infty}^{\psi}$ resp. les images par ces isomorphismes de $\widetilde{\mathrm{Jac}}_n^{\psi}$ et $\widetilde{\mathrm{Jac}}_{\infty}^{\psi}$. Alors $\widetilde{\mathcal{G}}_n^{\psi}$ (resp. $\widetilde{\mathcal{G}}_{\infty}^{\psi}$) est l'idéal de Fitting initial du $\mathbb{Z}_p[\psi][\Gamma_n]$-module $\widetilde{\mathcal{U}}_n^{\psi}/\widetilde{\mathrm{Jac}}_n^{\psi}$ (resp. du $\Lambda_{\psi}$-module $\widetilde{\mathcal{U}}_{\infty}^{\psi}/\widetilde{\mathrm{Jac}}_{\infty}^{\psi}$). En particulier, ces idéaux sont indépendants du choix des isomorphismes de départ.
\end{Lemme}

\begin{proof}[\normalfont{\bfseries{Preuve.}}]
C'est évident, puisque $$\widetilde{\mathcal{U}}_{\infty}^{\psi}/\widetilde{\mathrm{Jac}}_{\infty}^{\psi}\simeq\Lambda_{\psi}/\widetilde{\mathcal{G}}_{\infty}^{\psi}\ \mathrm{et}\   \widetilde{\mathcal{U}}_{n}^{\psi}/\widetilde{\mathrm{Jac}}_{n}^{\psi}\simeq\mathbb{Z}_p[\psi][\Gamma_n]/\widetilde{\mathcal{G}}_{n}^{\psi}.$$
\end{proof}

Alors, par le théorème 3.3.1., l'application $f_{\infty}$ définie en (2) induit une suite exacte de $\Lambda_{\psi}$-modules:
$$0\rightarrow \mathrm{tor}_{\Lambda_{\psi}}\widetilde{\mathfrak{X}}_{\infty}^{\psi} \rightarrow X_{\infty}^{\psi} \rightarrow\Lambda_{\psi}/(\mathcal{F}_{\psi})\rightarrow\Lambda_{\psi}/\widetilde{\mathcal{G}}_{\infty}^{\psi}\rightarrow 0.$$
\\De plus (voir appendice sur le "miroir"), la torsion du module de Bertrandias-Payan peut être vue comme un adjoint:
$$\mathrm{tor}_{\Lambda_{\psi}}\widetilde{\mathfrak{X}}_{\infty}^{\psi}\simeq\alpha(X'_{\infty}(-1))^{\psi}.$$
où $\alpha(.)$ est le foncteur "adjoint" d'Iwasawa. D'après les propriétés fonctorielles de l'adjoint, on a un pseudo-isomorphisme:
$$\alpha(X'_{\infty}(-1))^{\psi}\approx(X_{\infty}^{'\psi^{*}})^{\sharp}$$
où $(.)^{\sharp}$ désigne \textbf{l'involution d'Iwasawa} définie par $\sigma \mapsto \kappa(\sigma)\sigma^{-1}$ et $\approx$ signifie "pseudo-isomorphisme".
\\Donc le théorème 3.3.1. peut se réénoncer ainsi:

\begin{Theoreme} Pour tout $\psi$ impair, distinct de $\omega$,
on a une suite exacte de $\Lambda_{\psi}$-modules:
$$0\rightarrow \alpha(X'_{\infty}(-1))^{\psi} \rightarrow X_{\infty}^{\psi} \rightarrow \Lambda_{\psi}/(\mathcal{F}_{\psi})\rightarrow \Lambda_{\psi}/\widetilde{\mathcal{G}}_{\infty}^{\psi}\rightarrow 0.$$
En particulier, $\mathrm{char}_{\Lambda_{\psi}}\alpha(X'_{\infty}(-1))^{\psi}=\mathrm{char}_{\Lambda_{\psi}}\Lambda_{\psi}/\widetilde{\mathcal{G}}_{\infty}^{\psi}$.
\end{Theoreme}

Pour passer des idéaux caractéristiques aux idéaux de Fitting, on utilisera le lemme suivant (qui généralise le théorème 1 de [Cor] et nous a été suggéré par le rapporteur):

\begin{Lemme}
Soit $M$ un $\Lambda_{\psi}$-module de torsion, notons $M^{o}$ son sous-module fini maximal. Alors
$\mathrm{Fitt}_{\Lambda_{\psi}}(M)=(\mathrm{char}_{\Lambda_{\psi}}(M)).\mathrm{Fitt}_{\Lambda_{\psi}}(M^{o}).$
\end{Lemme}

\begin{proof}[\normalfont{\bfseries{Preuve.}}]
Si $M$ est fini, le lemme est trivial. Si ce n'est pas le cas, posons $M'=M/M^{o}$. Alors $M'$ n'a pas de sous-module fini non nul, donc est de dimension projective au plus égale à un sur $\Lambda_{\psi}$ (voir le début du lemme 3.1.1.). Dans cette situation, on sait ([H], chapitre 3) que les idéaux de Fitting se comportent de façon multiplicative dans la suite exacte $0\rightarrow M^{o}\rightarrow M \rightarrow M' \rightarrow 0$, i.e. $\mathrm{Fitt}_{\Lambda_{\psi}}(M)=\mathrm{Fitt}_{\Lambda_{\psi}}(M^{o})\mathrm{Fitt}_{\Lambda_{\psi}}(M')$. On applique alors le lemme 3.1.1. à $M'$. 
\end{proof}

\begin{Theoreme}
Pour tout caractère impair $\psi$ de $\mathrm{Gal}(K/\mathbb{Q})$ tel que $\psi\neq \omega$, $(\widetilde{\mathcal{G}}_{\infty}^{\psi})^{\sharp}$ est l'idéal de Fitting de $X^{'\psi^*}_{\infty}$ sur l'algèbre $\Lambda_{\psi^*}$. Si $\psi^*(p)\neq 1$, $(\mathcal{G}_{\infty}^{\psi})^{\sharp}$ est l'idéal Fitting de $X_{\infty}^{\psi^*}$.
\end{Theoreme} 

\begin{proof}[\normalfont{\bfseries{Preuve.}}]
On sait déjà par le théorème 3.4.1. que $(X^{'\psi^*}_{\infty})^{\sharp}$ et $\Lambda_{\psi}/(\widetilde{\mathcal{G}}_{\infty}^{\psi})$ ont même idéal caractéristique. Pour appliquer le lemme 3.4.1., il faut montrer que leurs sous-modules finis maximaux ont même idéal de Fitting sur $\Lambda_{\psi}$. Pour simplifier l'écriture, posons $Y_{\infty}=\widetilde{\mathcal{U}}_{\infty}^{\psi}/\widetilde{Jac}_{\infty}^{\psi}\simeq \Lambda_{\psi}/\widetilde{\mathcal{G}}_{\infty}^{\psi}$. Comme $Y_{\infty}$ et $(X^{'\psi^*}_{\infty})^{\sharp}$ ont même série caractéristique, la conjecture de Gross (qui est vraie pour les corps abéliens) entraîne que le sous-module fini maximal $Y_{\infty}^{o}$ est égal à $Y_{\infty}^{\Gamma^{p^n}}$ pour $n$>>0. 
\\Faisons la descente en supposant d'abord $\psi(p)\neq 1$. La suite exacte tautologique
$$0\rightarrow \widetilde{\mathrm{Jac}}^{\psi}_{\infty}\rightarrow \widetilde{\mathcal{U}}_{\infty}^{\psi} \rightarrow Y_{\infty} \rightarrow 0$$
donne par co-descente un diagramme commutatif:
\\
\xymatrix{
 & 0 \ar@{->}[r] & Y_{\infty}^{\Gamma^{p^n}} \ar@{->}[r] & (\widetilde{\mathrm{Jac}}_{\infty}^{\psi})_{\Gamma^{p^n}} \ar@{->}[r] \ar@{->}[d] & (\widetilde{\mathcal{U}}_{\infty}^{\psi})_{\Gamma^{p^n}} \ar@{->}[r] \ar@{^(->>}[d] & ... \\
 & & 0 \ar@{->}[r]  & \widetilde{\mathrm{Jac}}_{n}^{\psi} \ar@{->}[r]  & \widetilde{\mathcal{U}}_n^{\psi} \ar@{->}[r] & ...
 }
\ \\où l'isomorphisme vertical de droite provient de l'hypothèse $\psi(p)\neq 1$. Donc, pour $n$>>0, $Y_{\infty}^{o}$ s'identifie au noyau du morphisme naturel $(\widetilde{\mathrm{Jac}}_{\infty}^{\psi})_{\Gamma^{p^n}}\rightarrow \widetilde{\mathrm{Jac}}_n^{\psi}$, c'est-à-dire, d'après le corollaire 1, au noyau du morphisme naturel $\mathrm{fr}_{\Lambda_{\psi}}\mathfrak{X}_{\infty}^{\psi}\rightarrow \mathrm{fr}_{\mathbb{Z}_p} \mathfrak{X}_n^{\psi}$. Or la détermination de ce noyau résulte d'un calcul classique de la théorie d'Iwasawa, que nous rappelons pour la commodité du lecteur: d'après le théorème de structure, on a une suite exacte 
$$0\rightarrow \mathrm{fr}_{\Lambda_{\psi}}\mathfrak{X}_{\infty}^{\psi} \rightarrow L \rightarrow H \rightarrow 0,$$
où $L$ est $\Lambda_{\psi}$-libre et $H$ est fini ; par co-descente, on obtient une suite exacte
$$0\rightarrow H^{\Gamma^{p^n}} \rightarrow (\mathrm{fr}_{\Lambda_{\psi}} \mathfrak{X}_{\infty}^{\psi})_{\Gamma^{p^n}} \rightarrow L_{\Gamma^{p^n}} \rightarrow ...$$
qui montre, pour $n$>>0, que $H$ s'identifie à la $\mathbb{Z}_p$-torsion de $(\mathrm{fr}_{\Lambda_{\psi}}\mathfrak{X}_{\infty}^{\psi})_{\Gamma^{p^n}}$. Or, dans le diagramme commutatif naturel
\\
\xymatrix{
 & 0 \ar@{->}[r] & \mathrm{tor}_{\mathbb{Z}_p} \mathfrak{X}_n^{\psi} \ar@{->}[r] & \mathfrak{X}_n^{\psi} \ar@{->}[r]  & \mathrm{fr}_{\mathbb{Z}_p} \mathfrak{X}_n^{\psi} \ar@{->}[r] & 0 \\
 & & & (\mathfrak{X}_{\infty}^{\psi})_{\Gamma^{p^n}} \ar@{->}[r] \ar@{->}[u] & (\mathrm{fr}_{\Lambda_{\psi}}\mathfrak{X}_{\infty}^{\psi})_{\Gamma^{p^n}} \ar@{->}[u] & 
 }
\ \\la flèche verticale de gauche est un isomorphisme pour tout caractère $\psi$ non trivial, ce qui montre que $H$ est le noyau de la flèche verticale de droite, i.e. $Y_{\infty}^{o}\simeq H$. Mais il est bien connu que $H$ est canoniquement isomorphe au dual de Kummer du sous-module fini maximal $(X^{'\psi^*}_{\infty})^{o}$ (une référence implicite est [Iw2], une référence explicite [J], coroll 6.5):
$$H\simeq\mathrm{Hom}((X^{'\psi^*}_{\infty})^{o},\mu_{p^{\infty}})=\mathrm{Hom}((X^{'\psi^*}_{\infty})^{o}(-1),\mathbb{Q}_p/\mathbb{Z}_p).$$
Il est aussi connu (et cela se démontre facilement à partir de [MW], App. 4, propos. 1), que pour tout groupe fini $G$ dont le $p$-groupe de Sylow est cyclique, pour tout $\mathbb{Z}_p[G]$-module fini $M$, le $\mathbb{Z}_p[G]$-Fitting de $\mathrm{Hom}(M,\mathbb{Q}_p/\mathbb{Z}_p)$ s'obtient exactement à partir de celui de $M$ en inversant l'action de $G$. Comme les morphismes de liaison entre les modules considérés sont surjectifs, les idéaux de Fitting de ces modules sont compatibles avec la limite projective (voir e.g. [GK], thm. 2.1) et l'on obtient le résultat cherché, à savoir que $(X^{'\psi^*}_{\infty})^{\sharp}$ et $\Lambda_{\psi}/(\widetilde{\mathcal{G}}_{\infty}^{\psi})$ ont même idéal de Fitting sur $\Lambda_{\psi}$. La première partie du théorème 3.4.2. est ainsi démontrée.
\\Examinons maintenant le cas où $\psi(p)=1$. Comme on l'a déjà signalé (voir la remarque après le théo. 3.3.1.), $(X_{\infty}^{\psi})^{\Gamma^{p^n}}$ n'est plus fini et $\widetilde{\mathrm{Jac}}_n^{\psi}$ n'est plus isomorphe à $\mathrm{fr}_{\mathbb{Z}_p}\mathfrak{X}_n$. Cependant, le diagramme de co-descente
\\
\xymatrix{
 & 0 \ar@{->}[r] & Y_{\infty}^{\Gamma^{p^n}} \ar@{->}[r] & (\widetilde{\mathrm{Jac}}_{\infty}^{\psi})_{\Gamma^{p^n}} \ar@{->}[r] \ar@{->}[d] & (\widetilde{\mathcal{U}}_{\infty}^{\psi})_{\Gamma^{p^n}} \ar@{->}[r] \ar@{->}[d] & ... \\
 & & 0 \ar@{->}[r]  & \widetilde{\mathrm{Jac}}_{n}^{\psi} \ar@{->}[r]  & \widetilde{\mathcal{U}}_n^{\psi} \ar@{->}[r] & ...
 }
\ \\reste valable, et dans ce diagramme, l'homorphisme surjectif $\varphi : \ (\widetilde{\mathrm{Jac}}_{\infty}^{\psi})_{\Gamma^{p^n}}\simeq (\mathrm{fr}_{\Lambda_{\psi}}\mathfrak{X}_{\infty}^{\psi})_{\Gamma^{p^n}}\rightarrow \widetilde{\mathrm{Jac}}_n^{\psi}$ (sous-section 2.4) se factorise à travers $\mathrm{fr}_{\mathbb{Z}_p}\mathfrak{X}_n^{\psi}$. Comme précédemment, on a $H=\mathrm{tor}_{\mathbb{Z}_p}(\widetilde{\mathrm{Jac}}_{\infty}^{\psi})_{\Gamma^{p^n}}$ et la factorisation de $\varphi$ à travers $\mathrm{fr}_{\mathbb{Z}_p}\mathfrak{X}_n^{\psi}$ montre que $H=\mathrm{tor}_{\mathbb{Z}_p}(\mathrm{Ker} \varphi)$. Or une chasse dans le diagramme de co-descente montre immédiatement que $Y_{\infty}^{\Gamma^{p^n}}$ (qui est fini par Leopoldt) s'injecte dans $H$ et que le quotient $H/Y_{\infty}^{\Gamma^{p^n}}$ s'injecte dans $\mathrm{Ker}((\widetilde{\mathcal{U}}_{\infty})_{\Gamma^{p^n}}\rightarrow \widetilde{\mathcal{U}}_n)$. Reste seulement à déterminer ce dernier noyau. La suite exacte tautologique 
$$0\rightarrow \widetilde{\mathcal{U}}_{\infty}^{\psi}\rightarrow \mathfrak{X}_{\infty}^{\psi}\rightarrow X_{\infty}^{\psi}\rightarrow 0$$ donne par co-descente un diagramme commutatif:
\\
\xymatrix{
 0 \ar@{->}[r] & (X_{\infty}^{\psi})^{\Gamma^{p^n}} \ar@{->}[r] & (\widetilde{\mathcal{U}}_{\infty}^{\psi})_{\Gamma^{p^n}} \ar@{->}[r] \ar@{->}[d] & (\widetilde{\mathfrak{X}}_{\infty}^{\psi})_{\Gamma^{p^n}} \ar@{->}[r] \ar@{^(->>}[d] & (X_{\infty}^{\psi})_{\Gamma^{p^n}} \ar@{->}[r] \ar@{->}[d] & 0 \\
 & 0 \ar@{->}[r]  & \widetilde{\mathcal{U}}_{n}^{\psi} \ar@{->}[r]  & \widetilde{\mathfrak{X}}_n^{\psi} \ar@{->}[r] & A_n^{\psi} \ar@{->}[r] & 0
 }
\ \\qui montre que le noyau cherché coïncide avec $(X_{\infty}^{\psi})^{\Gamma^{p^n}}$. Comme $\psi$ est impair, on sait que $X_{\infty}^{\psi}$ n'a pas de sous-module fini non nul et en définitive, $H=Y_{\infty}^{\Gamma^{p^n}}$.
\\Dans le cas où $\psi^{*}(p)\neq 1$, on a par définition $\widetilde{\mathcal{G}}_{\infty}^{\psi}=\mathcal{G}_{\infty}^{\psi}$ et il est connu (voir par ex. [HI], pp. 380-381) que la surjection naturelle $X^{\psi^*}_{\infty}\rightarrow X^{'\psi^*}_{\infty}$ est alors un isomorphisme, d'où la seconde assertion du théorème.
\end{proof}
(Nous remercions Anthony Martin de nous avoir signalé une erreur dans une première version et de nous avoir suggéré la présente preuve.)

\begin{Corollaire}
$X^{'\psi^*}_{\infty}=(0) \Leftrightarrow \widetilde{\mathcal{G}}_{\infty}^{\psi}=\Lambda_{\psi} \Leftrightarrow $ l'application dans la suite exacte du théo. 3.4.1 réalise un isomorphisme $X_{\infty}^{\psi}\simeq \Lambda_{\psi}/(\mathcal{F}_{\psi})$.
\end{Corollaire}

\textbf{Remarque:} C'est l'analogue exact de la conjecture de Vandiver - sous différentes formulations - au-dessus de $\mathbb{Q}(\zeta_{p^{\infty}})$; rappelons que la conjecture de Greenberg prévoit seulement la finitude de $X^{'\psi^*}_{\infty}$.

\begin{proof}[\normalfont{\bfseries{Preuve.}}]
La première équivalence résulte immédiatement du théorème 3.4.2. La seconde équivalence résulte de la première et de la suite exacte du théorème 3.4.1.
\end{proof}

\begin{Corollaire} Soit $\psi$ impair, $\psi\neq \omega$.
\begin{enumerate}
	\item Si $\psi(p)\neq 1$, $(\widetilde{\mathcal{G}}_n^{\psi})^{\sharp}$ est l'idéal de Fitting de $X_{\infty}^{'\psi^*}/\omega_n$ sur l'algèbre $\mathbb{Z}_p[\psi^*][\Gamma_n]$. 
	\item Si en outre $\psi^*(p)\neq 1$, $(\mathcal{G}_n^{\psi})^{\sharp}$ est l'idéal de Fitting de $A_n^{\psi^*}$, $\forall n\geq 0$.
\end{enumerate}
\end{Corollaire}

\begin{proof}[\normalfont{\bfseries{Preuve.}}]
1. L'idéal de Fitting de $X^{'\psi^*}_{\infty}/\omega_n$ sur $\Lambda_{\psi^*}/\omega_n\simeq\mathbb{Z}_p[\psi^*][\Gamma_n]$ est l'image dans cette algèbre de $(\widetilde{\mathcal{G}}_{\infty}^{\psi})^{\sharp}$. Il reste à montrer la surjectivité de l'application naturelle $N: \ (\widetilde{\mathcal{G}}_{\infty}^{\psi})^{\sharp}\rightarrow (\widetilde{\mathcal{G}}_{n}^{\psi})^{\sharp}$ i.e., par le lemme de Nakayama, la surjectivité de  $N: \ (\widetilde{\mathcal{G}}_{\infty}^{\psi})^{\sharp}/p \rightarrow (\widetilde{\mathcal{G}}_{n}^{\psi})^{\sharp}/p$. Comme $\kappa(\sigma) \equiv 1\ \mathrm{mod}\ (p)$ pour tout $\sigma\in\Gamma^{p^n}$, $M^{\sharp}/p\simeq M/p$ pour tout $\Gamma^{p^n}$-module $M$, et la surjectivité cherchée équivaut à celle de  $N: \ \widetilde{\mathcal{G}}_{\infty}^{\psi}/p\rightarrow \widetilde{\mathcal{G}}_{n}^{\psi}/p$ i.e., par Nakayama, à la surjectivité de  $N: \ \widetilde{\mathcal{G}}_{\infty}^{\psi}\rightarrow \widetilde{\mathcal{G}}_{n}^{\psi}$, qui est claire par le paragraphe 2.4. La première assertion du corollaire est ainsi démontrée.
\\2. Par l'hypothèse de semi-simplicité, toutes les $p$-places sont totalement ramifiées dans $K_{\infty}/K$, donc les applications naturelles $X^{'\psi^*}_{\infty}/\omega_n\rightarrow A^{'\psi^*}_n$ sont surjectives. L'hypothèse $\psi^*(p)\neq 1$ entraîne qu'elles sont injectives d'après [LMN], proposition 4.1. Elle entraîne aussi, comme on l'a déjà vu dans le théorème 3.4.1., que $X_{\infty}^{\psi^*}\simeq X_{\infty}^{'\psi^*}$. 
\end{proof}

Bien entendu, l'intérêt du théorème 3.4.2. dépend de notre capacité à décrire plus ou moins explicitement l'idéal $\widetilde{\mathcal{G}}_{\infty}^{\psi}$. Le module des sommes de Jacobi $\overline{\mathrm{Jac}}_{\infty}^{\psi}$ étant isomorphe à $\mathcal{F}_{\psi}.\mathfrak{X}_{\infty}^{\psi}$ (lemme 3.2.2.), on est ramené à la détermination des éléments $\mathcal{F}_{\psi}.\rho$ pour $\rho\in\mathfrak{X}_{\infty}^{\psi}$, ce qui constitue en quelque sorte une loi de réciprocité explicite au niveau infini. Le principe du calcul sera le suivant: puisque $\psi\neq\omega$, la $\psi$-partie de $\mathcal{U}_{\infty}^{\psi}$ s'injecte dans $\mathfrak{X}_{\infty}^{\psi}$ (par le corps de classes, on peut l'identifier au sous-groupe d'inertie de $\mathfrak{X}_{\infty}^{\psi}$ en toutes les $p$-places). Comme $\mathcal{F}_{\psi}$ annule $X_{\infty}^{\psi}$, $\mathcal{F}_{\psi}.\rho\in\mathcal{U}_{\infty}^{\psi}$ pour tout $\rho\in \mathfrak{X}_{\infty}^{\psi}$, et l'on peut chercher à déterminer $\mathcal{F}_{\psi}.\rho$ par la théorie de Coleman, ce qui revient \textsl{in fine} à calculer certains caractères de Coates-Wiles. Le problème consiste en fait à établir pour les sommes de Gauss, dans la partie "moins" des unités semi-locales, un analogue du théorème classique d'Iwasawa (généralisé par Tsuji, voir [T2]) qui relie des unités cyclotomiques spéciales, dans la partie "plus" des unités semi-locales, aux éléments de Stickelberger. Par une approche "motivique" et en utilisant la loi de réciprocité de Bloch-Kato, on peut se ramener aux caractères de Soulé, c'est-à-dire à des caractères kummeriens construits à partir d'un système cohérent particulier de $p$-unités cyclotomiques. On obtient ainsi une description remarquablement simple et explicite des idéaux de Fitting de $X_{\infty}^{\psi^*}$ et des $A_n^{\psi^*}$ (théorème 5.3.2. ci-dessous).  

\section{Calculs explicites sur $\mathbb{Q}$}

Pour des raisons didactiques, on étudiera d'abord le cas particulier $F=\mathbb{Q}$, qui présente le double avantage, d'une part de ne pas obscurcir les arguments principaux par excès de détails techniques, d'autre part "d'expliquer" l'intervention des sommes de Gauss et de Jacobi. Mais à strictement parler, le seul résultat technique nécessaire aux calculs explicites est l'énoncé 2 du théorème 4.1.1. ci-après. Le lecteur intéressé seulement au cas général pourra se reporter directement au §5.
\ \\
Dans toute la suite, on désignera par $\omega$ le caractère de Teichmüller de $G=G_0=\mathrm{Gal}(\mathbb{Q}(\zeta_p)/\mathbb{Q})$. Tout caractère impair $\psi$ de $G$ s'écrira sous la forme $\psi=\omega^i$, $1\leq i \leq p-2$, $i$ impair. Si $M$ est un module sur lequel $G$ opère, on notera $M^{(i)}$ le sous-espace propre correspondant au caractère $\omega^{i}$.
\\Suivant la démarche générale exposée à la fin du §3.4, il nous faut choisir un homomorphisme "explicite" $$\mathcal{U}_{\infty}^{(i)}\rightarrow\Lambda_{\omega_i}\simeq\mathbb{Z}_p[[T]]e_i,$$
où $e_i$ est l'idempotent associé à $\omega^i$. 
\\Partons de l'homorphisme de Coleman $\mathcal{U}_{\infty}\rightarrow\mathbb{Z}_p[[G_{\infty}]].(1+T)$ défini par $$\epsilon\in\mathcal{U}_{\infty}\mapsto \mathrm{Col}(\epsilon)=(1-\varphi/p)\mathrm{log}\ f_{\epsilon}(T),$$ 
où $f_{\epsilon}(T)\in\mathbb{Z}_p((T))^{\times}$ est la série de Coleman associée à $\epsilon$ et $(\varphi g)(T)=g((1+T)^p-1)$. Rappelons que $G_{\infty}=\mathrm{Gal}(\mathbb{Q}(\zeta_{p^{\infty}})/\mathbb{Q})$ opère sur $\mathcal{U}_{\infty}$ naturellement, et sur $\mathbb{Z}_p[[T]]$ par $(\tau f)(T)=f((1+T)^{\kappa(\tau)}-1)$, où $\kappa$ désigne le caractère cyclotomique. Pour $\epsilon\in\mathcal{U}_{\infty}$, il existe alors un unique élément $\mathrm{Mel}(\epsilon)\in\mathbb{Z}_p[[G_{\infty}]]$ tel que $\mathrm{Col}(\epsilon)=\mathrm{Mel}(\epsilon).(1+T)$, ce qui permet de définir un $\mathbb{Z}_p[[G_{\infty}]]$-homomorphisme
$$\mathrm{Mel}:\ \mathcal{U}_{\infty}\rightarrow\mathbb{Z}_p[[G_{\infty}]]\simeq\mathbb{Z}_p[G][[T]].$$
D'après les considérations du §3.4, il est clair que l'idéal de Fitting $\mathcal{G}_{\infty}^{(i)}\subset\Lambda_{\omega^{i}}$ est formé des éléments $\mathrm{Mel}(\mathcal{F}_i .\rho)$, où $\mathcal{F}_{i}$ est la série caractéristique de $X_{\infty}^{(i)}$ et $\rho$ parcourt $\mathfrak{X}_{\infty}^{(i)}$. Il est connu que $X_{\infty}^{(1)}=(0)$, i.e. $\mathcal{F}_1$ est inversible.
\ \\
\\ \textbf{Remarque:} Les notations de la théorie de Coleman ne sont pas très bien fixées. On a adopté ici la terminologie de [PR], qui se justifie par le fait que Mel provient essentiellement d'une transformée de Mellin ([PR], p.91).
\\On va calculer $\mathrm{Mel}(\mathcal{F}_{i} \rho)$dans le cadre de la théorie d'Ihara sur les pro-$p$-revêtements étales de la droite projective privée de trois points.
\ \\
\\ \textbf{Rappels sur la théorie d'Ihara}:([Col2], [Ih])
\\ Le pro-$p$-revêtement abélien étale maximal de $\mathbb{P}_{\overline{\mathbb{Q}}}^{1}-\{0,1,\infty\}$ est donné par la tour des courbes de Fermat d'équations homogènes $X^{p^n}+Y^{p^n}+Z^{p^n}=0$, $n\geq1$, et le groupe de Galois absolu $\mathrm{G}_{\mathbb{Q}}$ opère naturellement sur le module de Tate $\mathbb{T}:=\displaystyle{\lim_{\longleftarrow}}\ T_p(J_n(\overline{\mathbb{Q}}))$ associé aux jacobiennes $J_n$ de ces courbes. La théorie d'Ihara permet d'en déduire des représentations galoisiennes décrites en termes de séries formelles. Citons trois résultats principaux:

\begin{Theoreme}
$\mathbb{T}$ est un module monogène sur l'anneau $$\mathcal{A}=\mathbb{Z}_p[[u,v]]\simeq\mathbb{Z}_p[[u,v,w]]/((1+u)(1+v)(1+w)-1).$$
\end{Theoreme}

Un générateur privilégié $\eta$ de $\mathbb{T}$ sur $\mathcal{A}$ provient du commutateur $xyx^{-1}y^{-1}$ dans le pro-$p$-groupe libre à deux générateurs $x$ et $y$, qui est le groupe fondamental du pro-$p$-revêtement étale maximal de $\mathbb{P}_{\overline{\mathbb{Q}}}^{1}-\{0,1,\infty\}$ (ce commutateur correspond au contour de Pochhammer dans $\mathbb{P}_{\mathbb{C}}^{1}-\{0,1,\infty\}$). En fixant une fois pour toutes un générateur de $T_p(\mathbb{G}_{m})=\mathbb{Z}_p(1)$, l'action de $\mathrm{G}_{\mathbb{Q}}$ sur $\mathbb{T}$ se traduit par un $1$-cocycle continu
$$\mathrm{Ih}:\ \mathrm{G}_{\mathbb{Q}}\rightarrow\mathcal{A}^{*},$$
$\sigma\mapsto\mathrm{Ih}_{\sigma}(u,v,w)$ tel que $\sigma.\eta=\mathrm{Ih}_{\sigma}(u,v,w).\eta$. Par restriction, on obtient une représentation continue de $\mathfrak{X}_{\infty}$ dans $\mathcal{A}^{*}$:

\begin{Theoreme}
On a une représentation continue $\mathrm{Ih}:\mathfrak{X}_{\infty}\rightarrow\mathcal{A}^{*}$, $\rho\mapsto\mathrm{Ih}_{\rho}(u,v,w)$ tel que
$$\mathrm{Ih}_{\rho}(u,v,w)=\mathrm{exp}\displaystyle{\sum_{m\geq3,\ m\ \mathrm{impair}.}}\frac{\beta_m(\rho)}{m!}(U^m+V^m+W^m),$$
$1+u=\mathrm{exp}U,\ 1+v=\mathrm{exp}V,\ 1+w=\mathrm{exp}W$. En particulier, $\mathrm{Ih}_{\rho}(u,v,w)=\mathrm{G}_{\rho}(u)\mathrm{G}_{\rho}(v)\mathrm{G}_{\rho}(w)$, $(1+u)(1+v)(1+w)=1$, avec
$$\mathrm{G}_{\rho}(u)=\mathrm{exp}\displaystyle{\sum_{m\geq3,\ m\ \mathrm{impair}.}}\frac{\beta_m(\rho)}{m!}U^m$$
et des expressions analogues pour $\mathrm{G}_{\rho}(v)$, $\mathrm{G}_{\rho}(w)$. Pour tout $m\geq3$, impair, $\beta_m\in\mathrm{Hom}_{G_{\infty}}(\mathfrak{X}_{\infty},\mathbb{Z}_p(m))$.
\end{Theoreme}

La série $\mathrm{Ih}_{\rho}$ (resp. $\mathrm{G}_{\rho}$) est universelle pour les sommes de Jacobi (resp pour les sommes de Gauss) en un sens bien précis, voir [Ih], thm.7 (resp. [Col2], thm. 5.3). La factorisation de $\mathrm{Ih}_{\rho}$ dans le théorème précédent est l'analogue de la factorisation d'une somme de Jacobi à deux paramètres en un produit de trois sommes de Gauss. De plus, alors que les valeurs spéciales de $\mathrm{G}_{\rho}$ sont reliées aux sommes de Gauss, les coefficients de $\mathrm{G}_{\rho}$ sont reliés aux unités cyclotomiques. 

\subsection{Rappels sur quelques caractères $p$-adiques.}

Rappelons d'abord la définition de quelques caractères qui interviennent classiquement en théorie d'Iwasawa cyclotomique (voir e.g. [IS]). Fixons une fois pour toutes un générateur $\nu=(\zeta_{p^n})_{n\geq 1}$ de $\mathbb{Z}_p(1)$.
\\\ \\
-- \textbf{les caractères de Coates-Wiles} $\phi^{CW}_m \in \mathrm{Hom}_{G_{\infty}}(\mathcal{U_{\infty}},\mathbb{Z}_p(m))$ apparaissent comme coefficients dans le développement en série de $\mathrm{log}\ f_{\epsilon}$:
	$$\forall\epsilon\in\mathcal{U}_{\infty}, \mathrm{log}\ f_{\epsilon}(u)= \displaystyle{\sum_{m\geq 1,\ m\ \mathrm{impair}.}}\frac{\phi_m(\epsilon)}{m!} U^m\ , \ 1+u=\mathrm{exp}\ U,$$
	et $\phi_m^{CW}(\epsilon)=\phi_m(\epsilon)\otimes\nu^{\otimes m}$.
\\ \ \\
-- \textbf{les caractères de Soulé} $\chi_m\in\mathrm{Hom}_{G_{\infty}}(\mathfrak{X}_{\infty},\mathbb{Z}_p(m))$, sont des caractères kummeriens sur des systèmes projectifs de $p$-unités cyclotomiques:
	\\$\forall n\geq 1, \ \forall m \geq 1$, impair, posons $$\epsilon_n^{(m)}=\displaystyle{\prod_{1\leq a\leq p^n,\ (a,p)=1}}(1-\zeta_{p^n}^a)^{a^{m+1}}$$
	et définissons $\chi_m$ par $\zeta_{p^n}^{\chi_m(\rho)}=\{(\epsilon_n^{(m)})^{\frac{1}{p^n}}\}^{\rho-1},\ \forall n\geq1\ \forall\rho\in\mathfrak{X}_{\infty}$. Alors:

\begin{Theoreme} Pour tout $m\geq3$, $m$ impair, on a:
\\1. $\forall\epsilon\in\mathcal{U}_{\infty},\ \beta_m(\epsilon)=L_p(m,\omega^{1-m})\phi_m^{CW}(\epsilon)$.
\\2. $\forall i\ \mathrm{impair},\ 3\leq i\leq p-2,\ \forall \epsilon\in \mathcal{U}_{\infty}^{(i)},
\ \forall m \equiv i \ \mathrm{mod} \ (p-1), \chi_m(\epsilon)=(1-p^{m-1}) L_p(m,\omega^{1-m})\phi_m^{CW}(\epsilon)$
\\3. $\forall i\ \mathrm{impair},\ 3\leq i\leq p-2,\ \forall\rho\in\mathfrak{X}_{\infty}^{(i)},\ \forall m \equiv i \ \mathrm{mod} \ (p-1), \chi_m(\rho)=(1-p^{m-1})\phi_m^{CW}(\mathcal{F}_i.\rho)$
\\4. $\forall i\ \mathrm{impair},\ 3\leq i\leq p-2,\ \forall\rho\in\mathfrak{X}_{\infty}^{(i)},\ \forall m \equiv i\ \mathrm{mod} \  (p-1), \beta_m(\rho)=\phi_m^{CW}(\mathcal{F}_i.\rho)=(1-p^{m-1})^{-1}\chi_m(\rho)$

\end{Theoreme}

\begin{proof}[\normalfont{\bfseries{Principe de la preuve.}}]
1. provient de [Ih], thm. 10. La relation 2. a été montrée par Coleman en utilisant sa loi de réciprocité explicite ([Col1]). Donnons une idée de la preuve de 3. suivant la méthode simplifiée de [Ich2]:
\\Comme $\mathcal{F}_{i}$ est "la" série caractéristique de $X_{\infty}^{(i)}$, on a déjà vu que $\mathcal{F}_{i}.\rho\in \mathcal{U}_{\infty}^{(i)}$. Puisque $\mathrm{Hom}_{G_{\infty}}(\mathfrak{X}_{\infty}^{(i)},\mathbb{Z}_p(m))=0$ pour $m$ $\not \equiv$ $i\ \mathrm{mod}\ (p-1)$, on peut se limiter au cas où $m\equiv i\ \mathrm{mod}\ (p-1)$; alors, pour $\rho=\epsilon\in\mathcal{U}_{\infty}^{(i)}$, les relations 2. et 3. coïncident par Mazur-Wiles, donc 3. découle de 2. si la restriction $$\mathrm{Hom}_{G_{\infty}}(\mathfrak{X}_{\infty}^{(i)},\mathbb{Z}_p(m))\rightarrow \mathrm{Hom}_{G_{\infty}}(\mathcal{U}_{\infty}^{(i)},\mathbb{Z}_p(m))$$
est injective. Or, toujours d'après Mazur-Wiles, cette injectivité a lieu si et seulement si $L_p(m,\omega^{1-i})\neq0$. Mais il n'y a qu'un nombre fini d'entiers $m$ pour lesquels $L_p(m,\omega^{1-i})=0$, et donc la continuité en $m$ des caractères $\chi_m$ et $(1-p^{m-1})\phi_m^{CW}$ entraîne la validité de la relation 3.
\\La relation 4., quant à elle, repose sur l'injectivité (inconditionnelle) d'un autre homomorphisme de restriction, à savoir
$$\mathrm{Hom}_{G_{\infty}}(\mathfrak{X}_{\infty}^{(i)},\mathcal{B})\rightarrow\mathrm{Hom}_{G_{\infty}}(\mathcal{U}_{\infty}^{(i)},\mathcal{B}),$$
où $\mathcal{B}=\{g\in\mathcal{A}^{*};\ g(0,0)\equiv 1 \mathrm{mod} \ p\}$ (voir e.g. [Ich2]).
\end{proof}

\subsection{Sommes de Gauss}

Le théorème précédent détermine clairement l'image $\mathrm{Col}(\mathcal{F}_i.\rho)$ en termes kummeriens. Plus précisément:

\begin{Lemme}
$\forall i \geq 3, \ \mathrm{impair},\ \forall \rho\in\mathfrak{X}_{\infty}^{(i)},\ \mathrm{Col}(\mathcal{F}_i.\rho)$ est la "série universelle de Gauss" $$H_{\rho}^{i}(T):=\displaystyle{\sum_{m\geq 3, m\equiv i\ (p-1)}}\frac{\chi_m(\rho)}{m!} X^m,$$
avec $\mathrm{exp}\ X=1+T$.
\end{Lemme}

\begin{proof}[\normalfont{\bfseries{Preuve.}}] 
A partir de la formule évidente $D(\varphi g)=p \varphi(Dg)$, où $D$ est l'opérateur différentiel $(1+T)\frac{d}{dT}$, on obtient immédiatement le développement en série
$$\mathrm{Col}(\epsilon)=\displaystyle{\sum_{m\geq1,\ m\ \mathrm{impair}}}\frac{(1-p^{m-1})\phi_m^{CW}(\epsilon)}{m!}X^m,$$
$\mathrm{exp}\ X=1+T$, $\forall\epsilon\in\mathcal{U}_{\infty}$. Pour $\epsilon=\mathcal{F}_{i}.\rho$, on en déduit, d'après 3. du théorème précédent, que $\mathrm{Col}(\mathcal{F}_{i}.\rho)=H_{\rho}^{i}(T)$.
\end{proof}

Il reste à déterminer $\mathrm{Mel}(\mathcal{F}_i.\rho)$, ce qu'on va faire en introduisant, pour tout $\rho\in\mathfrak{X}_{\infty}$, une mesure $\delta(\rho)\in\mathbb{Z}_p[[G_{\infty}]]$ définie par Ihara-Kaneko-Yukinari  ([IKY], prop.2):

\newtheorem{DefLemme}{Définition-Lemme}[subsection]

\begin{DefLemme}
Pour tout $\rho \in\mathfrak{X}_{\infty}$, pour tout $n\geq1$ et tout $\overline{a}\in(\mathbb{Z}/p^n\mathbb{Z})^{*}$, fixons $b_{n,\overline{a}}\in\mathbb{Z}_p$ tel que $\{(1-\zeta_{p^n}^{\overline{a}})^{1/p^n}\}^{\rho-1}=\zeta_{p^n}^{b_{n,\overline{a}}}$ et posons 
$$\delta_n(\rho)=\displaystyle{\sum_{a=1, (a,p)=1}^{p^n-1}}b_{n,\overline{a}}\ a^{-1}\sigma_{a}\in \mathbb{Z}_p[G_{\infty}],$$
où $\sigma_a$ est l'unique élément de $G_{\infty}$ tel que $\kappa(\sigma_a)=a$. Alors $\delta(\rho)=\displaystyle{\lim_{n \rightarrow \infty}}\ \delta_n\in\mathbb{Z}_p[[G_{\infty}]]$ existe, et 
$$\delta(\rho).(1+T)=\displaystyle{\sum_{m\geq 1}}\frac{\chi_m(\rho)}{m!} X^m\ , \ \mathrm{exp}\ X=1+T.$$
\end{DefLemme}

\begin{proof}[\normalfont{\bfseries{Principe de la preuve:}}] 
On remarque qu'en fait $\delta(\rho)= \displaystyle{\lim_{\longleftarrow}}\ \pi_n(\delta_n)$, où $\pi_n$ est la projection de $\mathbb{Z}_p[[G_{\infty}]]$ dans $\mathbb{Z}/p^n\mathbb{Z}[G_n]$. De plus, par définition, $$\chi_m(\rho)\equiv\displaystyle{\sum_{a=1\ (a,p)=1}^{p^n-1}} b_{n,\bar{a}}.a^{m-1} \ \mathrm{mod}\ p^n .$$ 
Or il est immédiat que: 
$$\delta_n(\rho).(1+T)=\displaystyle{\sum_{a=1, (a,p)=1}^{p^n-1}}b_{n,\overline{a}}\ a^{-1}(1+T)^a=\displaystyle{\sum_{m=0}^{\infty}}\{\displaystyle{\sum_{a=1, (a,p)=1}^{p^n-1}}b_{n,\overline{a}}.a^{m-1}\}\frac{X^m}{m!}.$$
Le reste du calcul est purement algébrique.
\end{proof}

Il en résulte immédiatement que pour tout $i$ impair, $3\leq i \leq p-2$, pour tout $\rho\in\mathfrak{X}_{\infty}^{(i)}$, $\mathrm{Mel}(\mathcal{F}_i.\rho)=\delta(\rho)$. D'où:

\begin{Proposition}
Pour tout $i$ impair, $3\leq i \leq p-2$, $\mathcal{G}_{\infty}^{(i)}=\{\delta(\rho),\ \rho\in\mathfrak{X}_{\infty}^{(i)}\}.$
\end{Proposition}

En appliquant l'involution d'Iwasawa $\sharp$, on obtient que $\mathfrak{d}(\rho):=\delta(\rho)^{\sharp}$ existe et est la limite des $\mathfrak{d}_n(\rho)=\displaystyle{\sum_{a=1,(a,p)=1}^{p^n-1}} b_{n,\bar{a}}.\sigma_a\in \mathbb{Z}_p[G_{\infty}]$. Le théorème 3.4.2 et son corollaire 4 donnent alors:

\begin{Corollaire} Pour tout $i$ impair, $3\leq i\leq p-2$, l'idéal de Fitting de $X_{\infty}^{(1-i)}$ (resp. de $A_n^{(1-i)}$) sur $\Lambda$ (resp. sur $\mathbb{Z}_p[\Gamma_n]$) est formé des $\mathfrak{d}(\rho)$ (resp. des projections des $\mathfrak{d}(\rho)$) pour $\rho$ parcourant $\mathfrak{X}_{\infty}^{(i)}$.
\end{Corollaire}
Ces résultats seront généralisés dans la section 5.3.

\section{Calculs explicites sur un corps totalement réel}

On se propose dans cette section de déterminer explicitement les idéaux de Fitting $\widetilde{\mathcal{G}}_{\infty}^{\psi}$ en se plaçant dans la situation suivante:
\\$F$ est un corps abélien totalement réel, $\Delta=\mathrm{Gal}(F/\mathbb{Q})$, $p$ premier avec $|\Delta|$ et $p$ ne se ramifie pas dans $F$.
\\$K=K_0=F(\zeta_p)$, $K_n=F(\zeta_{p^{n+1}})$, $K_{\infty}=\displaystyle{\bigcup_{n\geq0}}K_n$, $\Delta=\mathrm{Gal}(F/\mathbb{Q})$, $G=\mathrm{Gal}(K/F)$, $G_n=\mathrm{Gal}(K_n/F)$, $G_{\infty}=\mathrm{Gal}(K_{\infty}/F)$, $\Gamma_n=\mathrm{Gal}(K_n/K)$, $\Gamma=\mathrm{Gal}(K_{\infty}/K)$.
\\Un caractère $\psi$ de $\mathrm{Gal}(K/\mathbb{Q})$ impair s'écrit $\chi\omega^i$, où $\chi$ est un caractère (forcément pair) de $\Delta$ et $i\in\mathbb{Z}$ est impair.
\\Notons que l'hypothèse de non ramification de $p$ dans $F$ n'est pas vraiment restrictive puisque tout corps abélien totalement réel peut être plongé dans un corps de la forme $K_n^{+}$. Concernant l'hypothèse de semi-simplicité, voir les commentaires suivant le théorème 3.3.1. ainsi que les calculs "globaux" du §5.4.
\\ \ \\L'extension des méthodes du §4 au cas général n'est pas immédiate. La théorie d'Ihara étant spécifique à $\mathbb{Q}(\zeta_p)$, on perd l'interprétation des sommes de Gauss ou de Jacobi en termes de séries universelles. Mais l'étape essentielle, à savoir la relation de Coleman (théorème 4.1.1.), peut être généralisée en adoptant une approche "motivique" qui permettra de passer des caractères de Coates-Wiles (définis via la théorie de Coleman) aux caractères de Soulé (dont la définition est purement kummerienne).

\subsection{Généralités}

Rappelons d'abord les préliminaires de la \textbf{théorie de Coleman} sur $F$ (voir e.g. [T1]). Soit $\Delta_p$ le groupe de décomposition de $p$ dans $\Delta$ et $\sigma_p\in\Delta$ l'automorphisme de Frobenius en $p$. Posons
$$\widehat{\mathcal{O}}_{F}=\mathcal{O}_{F}\otimes\mathbb{Z}_p\simeq\displaystyle{\prod_{v|p}}\mathcal{O}_v,$$
où $\mathcal{O}_v$ est l'anneau des entiers du complété $F_v$. Comme $F/\mathbb{Q}$ est non ramifiée en $p$, le $\mathbb{Z}_p[\Delta]$-module additif $\widehat{\mathcal{O}}_{F}$ est isomorphe à $\mathbb{Z}_p[\Delta]$.
\\L'anneau des séries formelles $\widehat{\mathcal{O}}_{F}[[T]]$ admet une action de $\Delta$ via les coefficients, ainsi qu'une action de $G_{\infty}$ définie par:
$$(\tau g)(T)=g((1+T)^{\kappa(\tau)}-1),$$
où $\kappa$ est le caractère cyclotomique.
\\Soit $\phi$ l'endomorphisme continu de $\widehat{\mathcal{O}}_{F}[[T]]$ tel que
$$(\phi g)(T)=g^{\sigma_p}((1+T)^p-1).$$
L'homomorphisme de Coleman est le $\mathbb{Z}_p[\Delta][[G_{\infty}]]$-homomorphisme $\mathrm{Col }:\mathcal{U}_{\infty}\rightarrow\widehat{\mathcal{O}}_{F}[[T]]$ défini par $\mathrm{Col}(\epsilon)=(1-\phi/p)\mathrm{log}f_{\epsilon}(T)$, où $f_{\epsilon}(T)\in\widehat{\mathcal{O}}_{F}((T))^{\times}$ est la série de Coleman associée à $\epsilon=(\epsilon_n)_{n\geq0}$, vérifiant $f_{\epsilon}(\zeta_{p^{n+1}}-1)=(\epsilon_n)^{\sigma_p^n}$ $\forall n \geq 0$.
\\L'image de Col est contenue dans un module libre $\widehat{\mathcal{O}}_{F}[[G_{\infty}]].(1+T)$, de sorte qu'on peut définir $\mathrm{Mel}(\epsilon)\in\widehat{\mathcal{O}}_{F}[[G_{\infty}]]$ par la relation $\mathrm{Mel}(\epsilon).(1+T)=\mathrm{Col}(\epsilon)$. Le résultat fondamental de Coleman est l'existence d'une suite exacte de $\mathbb{Z}_p[\Delta][[G_{\infty}]]$-modules:
\\
\xymatrix{0 \ar[r] & \mathbb{Z}_p[\Delta/\Delta_p](1) \ar@{->}[r] & \mathcal{U}_{\infty} \ar@{->}[r]^-{\mathrm{Mel}} & \widehat{\mathcal{O}}_{F}[[G_{\infty}]] \ar@{->}[r]^-{\kappa} &  \displaystyle{\frac{\widehat{\mathcal{O}}_{F}}{(\sigma_p-1)\widehat{\mathcal{O}}_{F}}}(1) \ar@{->}[r] & 0}
\ \\ qui provient par induction (i.e. en tensorisant au-dessus de $\mathbb{Z}_p[\Delta_p][[G_{\infty}]]$ par $\mathbb{Z}_p[\Delta][[G_{\infty}]]$) d'une suite exacte de $\mathbb{Z}_p[\Delta_p][[G_{\infty}]]$-modules:
\\
\xymatrix{& 0 \ar[r] & \mathbb{Z}_p(1) \ar@{->}[r] & \mathcal{U}_{v,\infty} \ar@{->}[r] & \mathcal{O}_v[[G_{\infty}]] \ar@{->}[r]^-{\kappa} &  \mathbb{Z}_p(1) \ar@{->}[r] & 0}
\ \\ (où $\mathcal{U}_{v, \infty}$ désigne la limite projective des unités locales principales).
\\On notera que pour $\psi\neq \omega$, $\mathbb{Z}_p[\Delta/\Delta_p](1)^{\psi}$ n'est autre que le module $\mathfrak{T}_{\infty}^{\psi}$ introduit à la fin du § 3.2. Pour un caractère $\psi$ impair comme au début du §5, $\Lambda_{\psi}:=\widehat{\mathcal{O}}_F[[G_{\infty}]]^{\psi}\simeq\mathbb{Z}_p[\chi][[T]]e_i$, $e_i$ étant l'idempotent associé à $\omega^i$. D'après la suite exacte de Coleman et la stratégie exposée à la fin du §3, il est clair que l'image par $\mathrm{Mel}$ de $\mathcal{F}_{\psi}.\mathfrak{X}_{\infty}^{\psi}\subset\mathcal{U}_{\infty}^{\psi}$ n'est autre que $\widetilde{\mathcal{G}}_{\infty}^{\psi}\subset\Lambda_{\psi}$. Pour calculer les éléments $\mathrm{Mel}(\mathcal{F}_{\psi}.\rho)$, $\rho\in\mathfrak{X}_{\infty}^{\psi}$, on dispose encore des caractères de Coates-Wiles et de Soulé dans le cas général, même s'ils sont plus compliqués à définir que dans le cas particulier du paragraphe précédent:
\\ \ \\
---Pour un corps $F$ totalement réel, non ramifié en $p$, \textbf{les caractères de Coates-Wiles locaux} 
	$$\phi_{m,v}^{CW}:\ \mathcal{U}_{v,\infty}\rightarrow F_v(m):=F_v\otimes\mathbb{Q}_p(m)$$  
relatifs à $F_v$ ($v|p$) sont donnés à partir du développement en série 
$$\mathrm{log}\ f_{\epsilon}(u)=\displaystyle{\sum_{m\geq 1}}\frac{\phi_{m,v}(\epsilon)}{m!}U^m,$$
où $\epsilon\in\mathcal{U}_{v,\infty}$, $f_{\epsilon}$ est la série de Coleman associée, $1+u=\mathrm{exp}U$. En choisissant un générateur $\nu$ de $\mathbb{Z}_p(1)$, on pose, pour tout $m\in\mathbb{Z}$ : $\phi_{m,v}^{CW}(\epsilon)=\phi_{m,v}(\epsilon)\otimes\nu^{\otimes m}$. Alors $\phi_{m,v}^{CW}\in\mathrm{Hom}_{G_{\infty}}(\mathcal{U}_{v,\infty},F_v(m))\simeq\mathrm{H}^1(F_v,F_v(m))$ ([BK], §2). En induisant de $\Delta_p$ à $\Delta$, on peut définir un \textbf{caractère de Coates-Wiles semi-local} relatif à $F$,
$$\phi_m^{CW}:\mathcal{U}_{\infty}\rightarrow F_p(m),$$
avec $F_p:=F\otimes \mathbb{Q}_p\simeq \oplus_{v|p}F_v$, qui appartient à $\mathrm{Hom}_{G_{\infty}}(\mathcal{U}_{\infty},F_p(m))\simeq \oplus_{v|p} \mathrm{H}^1(F_v,F_v(m)).$
\\ \ \\
---Les \textbf{caractères de Soulé} $\chi_m^{(e)}$ relatifs à un corps de nombres \textbf{quelconque} $F$ sont des caractères kummeriens attachés à un système projectif $e$ de $p$-unités: si $e=(e_n)_{n \geq 0}\in \overline{U}_{\infty}'= \displaystyle{\lim_{\longleftarrow}}(U'_n\otimes\mathbb{Z}_p)$, où $U'_n$ est le groupe des $p$-unités de $K_n=F(\zeta_{p^n})$, on pose 
$$\epsilon_n^{(m)}(e)=\displaystyle{\prod_{\tau}}e_n^{\tau.<\kappa(\tau)^{m-1}>_n},$$
où $\tau$ parcourt $\mathrm{Gal}(K_n/F)$ et $<\alpha>_n$ désigne l'unique entier de l'intervalle $[0,p^n]$ qui est congru à $\alpha$ mod $p^n$. On définit alors le caractère $\chi_m^{(e)}\in\mathrm{Hom}_{G_{\infty}}(\mathfrak{X}_{\infty},\mathbb{Z}_p(m))$ par ([Sou], [IS]):
$$(\zeta_{p^n})^{\chi_m^{(e)}(\rho)}=\{(\epsilon_n^{(m)}(e))^{1/p^n}\}^{\rho-1},\ \ \forall\rho\in\mathfrak{X}_{\infty},\ \forall n\geq1.$$
Une interprétation peut-être plus fonctorielle consiste à utiliser l'isomorphisme ([KNF], lemma 2.2)
$$\mathrm{Hom}_{G_{\infty}}(\mathfrak{X}_{\infty},\mathbb{Z}_p(m))\simeq\mathrm{H}^1_{\acute{e}t}(\mathcal{O}_F[1/p],\mathbb{Z}_p(m))\mathrm{\ mod\ torsion}$$
et pour $m \neq 0,1$, la suite exacte ([KNF], thm 3.2)
$$0\rightarrow\overline{U}_{\infty}'(m-1)_{G_{\infty}}\rightarrow\mathrm{H}^1_{\acute{e}t}(\mathcal{O}_F[1/p],\mathbb{Z}_p(m))\rightarrow X'_{\infty}(m-1)^{G_{\infty}}\rightarrow 0$$
(toutes les constructions sont canoniques). On peut alors définir le caractère $\chi_m^{(e)}$ comme l'image de $(e_n\otimes\zeta_{p^n}^{\otimes m-1})_n\in(\overline{U}_{\infty}'(m-1))_{G_{\infty}}$ dans $\mathrm{H}^1_{\acute{e}t}(\mathcal{O}_F[1/p],\mathbb{Z}_p(m))$ mod torsion.
\\Si $F$ est abélien, le caractère de Soulé proprement dit, noté $\chi_m=\chi_m(F)$, correspond à un choix particulier de $e=(e_n)_{n\geq 1}$: on prend pour $e_n$ la $p$-unité cyclotomique $\mathrm{N}_{\mathbb{Q}(\zeta_{fp^n})/K_n}(1-\zeta_{fp^n})$, où $f$ est le conducteur de $K$. Dans la suite, il sera même commode de considérer plutôt l'image $\chi_m^{S}$ (S pour Soulé) de $\chi_m$ dans $\oplus_{v|p}\mathrm{H}^1(F_v,\mathbb{Z}_p(m))$ mod torsion, via le morphisme de localisation
$$\mathrm{loc}_m:\ \mathrm{H}^1_{\acute{e}t}(\mathcal{O}_F[1/p],\mathbb{Z}_p(m))\rightarrow\oplus_{v|p}\mathrm{H}^1(F_v,\mathbb{Z}_p(m)).$$
\textbf{Remarques:}
\\1. Pour $m\geq 1$, la $\mathbb{Z}_p$-torsion de $\mathrm{H}^1_{\acute{e}t}(\mathcal{O}_F[1/p],\mathbb{Z}_p(m))$ n'est autre que  ([KNF], lemma 2.2) $\mathrm{H}^0_{\acute{e}t}(\mathcal{O}_F[1/p],\mathbb{Q}_p/\mathbb{Z}_p(m))$.
\\2. Si $F$ est abélien totalement réel, pour tout caractère $\chi$ de $\Delta=\mathrm{Gal}(F/\mathbb{Q})$, la $\chi$-composante $\mathrm{loc}_m^{\chi}$ de $\mathrm{loc}_m$ est injective si et seulement si la valeur spéciale $L_p(F,\chi\omega^{1-m},m)$ est non nulle ([KN] thm. 3.1), ce qui se produit pour presque tout $m$.	
\\3. Certains énoncés de [KN] et [KNF] (mais qui ne sont pas utilisés ici) comportent des facteurs euleriens erronés, qui sont corrigés dans [BenN]

Dans le cas particulier d'un corps cyclotomique $F=\mathbb{Q}(\zeta_N)$, le caractère de Soulé $\chi_m$ coïncide mod torsion avec l'\textbf{élément cyclotomique de Deligne-Soulé} $c_m(\zeta_N) \in\mathrm{H}^1_{\acute{e}t}(\mathbb{Z}[\zeta_N][1/p],\mathbb{Z}_p(m))$
$$c_m(\zeta_N):=\displaystyle{\lim_{\longleftarrow}}\ \mathrm{Cores}_{ \mathbb{Q}(\zeta_{Np^n})/\mathbb{Q}(\zeta_N)}(1-\zeta_{Np^n})\otimes(\zeta_{p^n})^{\otimes m-1}$$
(voir e.g [HK], def 3.1.2, où notre $m$ est noté $r$), qui intervient dans la théorie des valeurs spéciales des fonctions $L$.
\\Notre objectif va maintenant être de trouver une relation entre $\phi_m^{CW}$ et $\chi_m^{S}$ généralisant la formule de Coleman dans le théorème 4.1.1. Le calcul reposera sur la loi de réciprocité de Bloch-Kato et la compatibilité entre régulateurs étale et syntomique.
\\ \ \\\textbf{La loi de réciprocité de Bloch-Kato} ([BK], [PR]) pour le motif $\mathbb{Z}(m)$ peut être considérée comme une version "supérieure" des lois de réciprocité classiques pour le motif $\mathbb{Z}(1)$ ou le groupe multiplicatif $\mathbb{G}_m$. Soit $F$ totalement réel, non ramifié en $p$. Pour toute place $v|p$, pour tout $m\geq 1$, le cobord $\partial_v^m:\ F_v=\mathrm{H}^0(F_v,\mathrm{B}_{cris}^+)\rightarrow\mathrm{H}^1(F_v,\mathbb{Q}_p(m))$ associé à la suite exacte canonique
\\
\xymatrix{& & 0 \ar[r] & \mathbb{Q}_p(m) \ar@{->}[r] & \mathrm{J}_{\mathbb{Q}}^{[m]} \ar@{->}[rr]^-{1-p^{-m}\sigma_p}& & \mathrm{B}_{cris}^+ \ar@{->}[r] & 0}
\ \\
(pour les définitions manquantes, voir [BK, p.339]) vérifie:
$$\partial_v^m(a_v)=\frac{1}{(m-1)!}T(a_v.\phi_{m,v}^{CW}),$$
où $T:\ \mathrm{H}^1(F_v,F_v(m))\rightarrow\mathrm{H}^1(F_v,\mathbb{Q}_p(m))$ se déduit de la trace de $F_v/\mathbb{Q}_p$ ([BK], thm. 2-1).

\subsection{Calculs cyclotomiques.}

La tactique, dont l'idée figure déjà dans l'appendice de Kurihara à [Gr1], sous un langage un peu différent, va consister à appliquer la réciprocité de Bloch-Kato à un élément spécial convenablement choisi. Dans cette sous-section, on travaillera dans des corps cyclotomiques pour lesquels il convient de donner quelques notations supplémentaires. Fixons un caractère $\chi$ (forcément pair) de $\Delta=\mathrm{Gal}(F/\mathbb{Q})$, de conducteur $N=N_{\chi}$. Les caractères de Coates-Wiles et de Soulé seront relatifs à $\mathbb{Q}(\zeta_N)^+$ et notés $\phi_m^{CW}(N)$, $\chi_m(N)$, $\chi_m^{S}(N)$. Dans l'algèbre du groupe $\mathrm{Gal}(\mathbb{Q}(\zeta_N)/\mathbb{Q})\simeq(\mathbb{Z}/N \mathbb{Z})^{\times}$, on désignera par $e_{\chi}^N$ l'idempotent associé à $\chi$ et 
$\epsilon_{\chi}^N =\displaystyle{\sum_{\tau\in(\mathbb{Z}/N \mathbb{Z})^{\times}}}\chi(\tau)^{-1}\tau$ son numérateur
 (pour appliquer $e_{\chi}^N$ ou $\epsilon_{\chi}^N$ il faut évidemment élargir les scalaires, mais pour alléger le texte, on évitera autant que possible de faire figurer les scalaires dans les notations). Les éléments $\chi_m^N:=c_m(\zeta_N)$ sont compatibles pour la norme et définissent des systèmes d'Euler, et sont donc reliés aux fonctions $L_p$. En les modifiant par des facteurs euleriens, Kato a construit des éléments $\widetilde{c}_m(\zeta_N)$ tels que $c_m(\zeta_N)=(1-p^{m-1}\sigma_p^{-1})\widetilde{c}_m(\zeta_N)$, et ces $\widetilde{c}_m(\zeta_N)$ proviennent d'éléments spéciaux de la K-théorie. Plus précisément, pour $m \geq 3$, il existe des éléments spéciaux $b_m(\zeta_N)\in K_{2m-1}(\mathbb{Q}(\zeta_N))\otimes\mathbb{Q}$ (appelés \textbf{éléments de Beilinson}) tel que le \textbf{régulateur étale} de Soulé 
$$r^m_{\acute{e}t}:\ K_{2m-1}(\mathbb{Q}(\zeta_N))\otimes\mathbb{Q}\rightarrow\mathrm{H}^1_{\acute{e}t}(\mathbb{Z}[\zeta_N][1/p],\mathbb{Q}_p(m))\rightarrow V_{\acute{e}t}^m$$
(où la première flèche est une classe de Chern $p$-adique, la deuxième est $\mathrm{loc}_m$ et $V_m^{\acute{e}t}:=\oplus_{v|p}\mathrm{H}^1_{\acute{e}t}(\mathbb{Q}(\zeta_N)_v,\mathbb{Q}_p(m))$) envoie $b_m(\zeta_N)$ sur $\widetilde{c}_m(\zeta_N)/N^{m-1}(m-1)!$ (Deligne-Beilinson-Huber-Wildeshaus, voir e.g. [HK], thm. 5-2-2).
\\Les valeurs spéciales des fonctions $L_p$, quant à elles, interviennent via les \textbf{régulateurs syntomiques} de Gros-Somekawa (voir [Gr1], ou encore [KN], §1). Pour toute place $v|p$, on a un homomorphisme composé
\\
\xymatrix{&r_{syn,v}^m: K_{2m-1}(\mathbb{Q}(\zeta_N)) \ar@{->}[r] & \mathrm{H}^1(\mathbb{Z}_p[\zeta_N],S_{\infty}(m)_{\mathbb{Q}_p}) \ar@{->}[r]^-{\alpha_v}& \mathbb{Q}(\zeta_N)_v}
\ \\
où la première flèche est une classe de Chern syntomique à valeurs dans un groupe de cohomologie syntomique; si $p\not{|}N$, ce groupe est égal à $\mathbb{Q}(\zeta_N)_v$ et $\alpha_v=(1-p^{-m}\sigma_p)$.
\\Alors $b_m(\zeta_N)$ s'envoie via $r_{syn,v}^m$ sur $(1-p^{-m}\sigma_p)Li_m(\zeta_N)$ (voir e.g. [Gr2], thm. 2.22, ainsi que la remarque 1.6 de [KN], p.426), où $Li_m$ est le polylogarithme $p$-adique de Coleman, dont le lien avec les fonctions $L_p$ est donné, pour tout caractère de Dirichlet pair $\theta\neq 1$, de conducteur $N_{\theta}$, par la formule ([Col3], p.172):
$$(1-\theta(p)^{-1}p^{-m})^{-1}.L_p(\theta^{-1}\omega^{1-m},m)=\tau(\theta)N_{\theta}^{-1}\displaystyle{\sum_{b=1}^{N_{\theta}}}\theta(b)Li_m(\zeta_{N_{\theta}}^b)$$
où $\tau(\theta)=\displaystyle{\sum_{b=1}^{N_{\theta}}}\theta(b)\zeta_{N_{\theta}}^b$ est la somme de Gauss classique pour $\theta$. Si $p$ ne divise pas $N_{\theta}$, $\tau(\theta)$ est une unité $p$-adique, et l'on posera $\eta_{\theta}=\tau(\theta)^{-1}N_{\theta}$ (indépendant de $m$).
Pour $\theta=1$, la formule précédente doit être remplacée par une formule à la limite au voisinage de $1$ ([Col3], p.173). Dans la suite, pour alléger le texte, on ne fera pas les calculs pour $\theta=1$ (qui sont parallèles à ceux du cas $\theta\neq 1$); voir aussi le §4.
\\ \ \\
La compatibilité entre les régulateurs syntomique et étale (théorème de Gros-Niziol, voir e.g. [KN], thm. 3.7) se résume en un triangle commutatif:
\\ \ \\
\xymatrix{ &  & V^{syn}_m :=\oplus_{v|p}\alpha_v(\mathrm{H}^1(\mathbb{Z}_p[\zeta_N],S_{\infty}(m)_{\mathbb{Q}_p})) \ar@{->}[dd]^{\partial^m}  
\\ K_{2m-1}(\mathbb{Q}(\zeta_N))\otimes\mathbb{Q}  \ar@{->}[urr]_-{r^m_{syn}}  \ar@{->}[drr]_{r^m_{\acute{e}t}}& & 
\\& & V^{\acute{e}t}_m:=\oplus_{v|p}\mathrm{H}^1_{\acute{e}t}(\mathbb{Q}(\zeta_N)_v,\mathbb{Q}_p(m))}
\ \\
les flèches $r_{syn}^m$ et $\partial^m$ étant définies à partir de $r^m_{syn,v}$ et $\partial_v^m$ par induction de $H_p$ à $H$, où $H$ désigne momentanément $\mathrm{Gal}(\mathbb{Q}(\zeta_N)/\mathbb{Q})$ et $H_p$ son sous-groupe de décomposition en une place $v|p$. Notons que $V^{syn}_m=\oplus_{v|p}\mathbb{Q}(\zeta_N)_v$.
\\La tactique va consister à calculer les images de $\epsilon_{\chi^{-1}}^N(b_m(\zeta_N))$ dans le triangle précédent:

\ \\
-- \textbf{Par le régulateur étale:}
	\\ D'après les rappels, 
	$$r^m_{\acute{e}t}(\epsilon_{\chi^{-1}}^N b_m(\zeta_N))=\frac{1}{N^{m-1}(m-1)!} \ \epsilon_{\chi^{-1}}^N(\widetilde{c}_m(\zeta_N))$$
i.e. 
	$$r^m_{\acute{e}t}(\epsilon_{\chi^{-1}}^N b_m(\zeta_N))=\frac{(1-p^{m-1}\chi(p))^{-1}}{N^{m-1}(m-1)!}\  \epsilon_{\chi^{-1}}^N\chi_m^{S}(N).$$
Dans la suite, on aura besoin de l'action explicite de cet opérateur sur la $\psi$-partie $\mathcal{U}_{\infty}^{\psi}(N)$. Pour alléger les notations (uniquement dans ce calcul), on écrira $\chi_m^S$ pour $\chi_m^S(N)$. Pour tout $u\in \mathcal{U}_{\infty}^{\psi}(N)$, pour tout $\sigma\in H$, on a, en tenant compte de l'action galoisienne habituelle sur les homomorphismes:
$$(\sigma \chi_m^S)(u)=\sigma(\chi_m^S(\sigma^{-1}u))=\chi(\sigma)^{-1}\omega(\sigma)^{-i}\sigma(\chi_m^S(u))=\chi(\sigma)^{-1}\chi_m^S(u)$$
puisque $H$ opère trivialement sur $\mu_{p^{\infty}}$ par définition. Il s'ensuit que 
\begin{equation}
(\epsilon_{\chi^{-1}}^N \chi_m^S)(u)=\displaystyle{\sum_{\sigma\in H}}\chi(\sigma)(\sigma \chi_m^S)(u)=\displaystyle{\sum_{\sigma\in H}}\chi_m^S(u)=|H|.\chi_m^S(u).
\end{equation}

\ \\
-- \textbf{Par le régulateur syntomique:}
\\Fixons une place $v|p$ et notons $\chi_p$ la restriction de $\chi$ à $H_p$, $\epsilon_{\chi_p^{-1}}$ le numérateur de l'idempotent relatif à $\chi_p^{-1}$ et $H_p$. D'après les rappels, le régulateur $r_{syn,v}^m$ envoie $b_m(\zeta_N)$ sur $(1-p^{-m}\sigma_p)\epsilon_{\chi_p^{-1}}Li_m(\zeta_N)=(1-p^{-m}\chi(p)^{-1})\displaystyle{\sum_{t\in H_p}}\chi(t)t(Li_m(\zeta_N)).$ Comme $\epsilon_{\chi^{-1}}^N=(\displaystyle\sum_{H/H_p}\chi(s)s)(\displaystyle\sum_{H_p}\chi(t)t)$, où $s$ parcourt un système complet de représentants de $H/H_p$, le régulateur syntomique $r^m_{syn}$ va envoyer $\epsilon_{\chi^{-1}}^{N} b_m(\zeta_N)$ sur $(1-p^{-m}\chi(p)^{-1})\displaystyle{\sum_{s,t}}\chi(st)st(Li_m(\zeta_N))=(1-p^{-m}\chi(p)^{-1})\displaystyle{\sum_{b=1}^N}\chi(b)Li_m(\zeta_N^b)$ puisque le polylogarithme est Galois-équivariant. En appliquant la formule de Coleman, on obtient finalement (cp. [Gr2], propos. 3.1):

$$r_{syn}^m(\epsilon_{\chi^{-1}}^N(b_m(\zeta_N)))=\eta_{\chi}L_p(\chi^{-1}\omega^{1-m},m).$$

On rappelle que $\eta_{\chi}=\tau(\chi)^{-1}N$ est une unité $p$-adique (indépendante de $m$).

\ \\
-- \textbf{Par Bloch-Kato:} 
\\Comme précédemment, fixons une place $v|p$ et posons 
$$a_v=r_{syn,v}^m(\epsilon_{\chi_p^{-1}}(b_m(\zeta_N))).$$
Pour $u_v\in\mathcal{U}_{v,\infty}(N)$, la formule de Bloch-Kato s'écrit $(m-1)!\partial_v^m(a_v)(u_v)=\displaystyle{\sum_{t\in H_p}}t(a_v)t(\phi_{m,v}(u_v))=\displaystyle{\sum_{t\in H_p}}\chi(t)^{-1}a_v t(\phi_{m,v}(u_v))=a_v\epsilon_{\chi_p}(\phi_{m,v}(u_v))$,
où l'on a noté $\phi_{m,v}$ pour $\phi_{m,v}^{CW}(N)$. La quantité $\epsilon_{\chi_p}(\phi_{m,v}(u_v))$ appartient à la $\chi_p$-partie de $\mathbb{Q}(\zeta_N)_v (\chi_p)=\mathbb{Q}_p(\zeta_N)(\chi_p)$. On peut encore la simplifier en prenant $u\in\mathcal{U}_{\infty}(N)^{\psi}$, car alors $\epsilon_{\chi_p}(\phi_{m,v}(u_v))=\phi_{m,v}(\epsilon_{\chi_p}(u_v))=|H_p|.\phi_{m,v}(u_v)$ (comme dans le calcul étale). En multipliant par $\sum_{s\in H/H_p} \chi(s) s$ i.e. par induction, on obtient, compte tenu des calculs syntomiques:
\begin{equation}
(\sum_{s\in H/H_p} \chi(s) s)(a_v \epsilon_{\chi_p}\phi_{m,v}(u_v))=|H_p| \eta_{\chi} \phi_m^{CW}(N)(u) L_p(\chi^{-1}\omega^{1-m},m) \end{equation}
Le caractère miroir de $\chi^{-1}\omega^{1-m}$ est $\chi \omega^m$, qui est égal à $\psi$ si $m\equiv i$ mod $(p-1)$. En appliquant la définition par Iwasawa de la fonction $L_p$ ([T1], thm. 4.3.) et la Conjecture Principale (thm. de Mazur-Wiles), on peut alors choisir une série caractéristique $\mathcal{F}_{\psi}$ de $X_{\infty}^{\psi}$ telle que
$$\eta_{\chi}|H|^{-1}|H_p|\phi_m^{CW}(N)(u)L_p(\chi^{-1}\omega^{1-m})=\phi_m^{CW}(N)(\mathcal{F}_{\psi}.u).$$
En tenant compte de la commutativité du triangle de régulateurs, on obtient, pour $u$ élément de $\mathcal{U}_{\infty}(N)^{\psi}$:
\begin{equation}
\forall m \equiv i \ \mathrm{mod}\ (p-1),\ \chi_m^S(N)(u)=N^{m-1}(1-p^{m-1}\chi(p))\phi_m^{CW}(N)(\mathcal{F}_{\psi}.u)
\end{equation}
NB: pour $N=1$, on retrouve bien la formule de Coleman (théorème 4.1.1. 2.).

\subsection{Retour à $F$}

On revient aux hypothèses et notations générales du début de la partie 5, où $\chi$ est un caractère (forcément pair) de $\Delta$, de conducteur $N_{\chi}$. La relation de Coleman (théorème 4.1.1. (3)) se généralise comme suit:

\begin{Theoreme}
Soit $\psi=\chi \omega^{i},\ 1\leq i\leq p-2$, impair. On peut choisir une série caractéristique $\mathcal{F}_{\psi}$ de $X_{\infty}^{\psi}$telle que pour tout $m\equiv i \ \mathrm{mod}\ (p-1)$, pour tout $\rho\in\mathfrak{X}_{\infty}^{\psi}$, on a:
$$\chi_m^S(\rho)=N_{\chi}^{m-1}(1-p^{m-1}\chi(p))\phi_m^{CW}(\mathcal{F}_{\psi}.\rho)$$
(ici, les notations $\chi_m^{S}$, $\phi_m^{CW}$, $\mathcal{U}_{\infty}^{\psi}$ sont relatives au corps de base $F$).
\end{Theoreme}

\begin{proof}[\normalfont{\bfseries{Preuve.}}]
Il s'agit de montrer que la formule (5) des calculs cyclotomiques précédents s'étend au corps de base $F$. Remarquons de prime abord que pour tout $F$ abélien, $\mathcal{U}_{\infty}^{\psi}$ est de $\Lambda_{\psi}$-rang un, $\mathrm{Hom}_{G_{\infty}}(\mathcal{U}_{\infty}^{\psi},\mathbb{Z}_p(m))$ est de $\mathbb{Z}_p$-rang un, et par suite:
$$\forall \epsilon \in \mathcal{U}_{\infty}^{\psi},\ \chi_m^S(\epsilon)=\lambda \phi_m^{CW}(\mathcal{F}_{\psi}.\epsilon),$$
où $\lambda\in\mathbb{Q}_p(\chi)$ et les notations se réfèrent à $F$. Pour déterminer $\lambda$, il suffira de faire un calcul particulier sur $\epsilon$. 
\\Notons $N$ le conducteur de $\chi$, $M$ celui de $F$ et prenons d'abord $F=\mathbb{Q}(\zeta_M)^{+}$. Il s'agit en fait de refaire les calculs cyclotomiques précédents en considérant $\chi$ comme un caractère mod $M$:
\\--à cause de la relation de distribution
$$M^{m-1}\sum_{a=1}^M \chi(a) Li_m(\zeta_M^a)=\prod_{l|M, (l,N)=1}(1-\chi(l)l^{m-1})N^{m-1}\sum_{b=1}^N \chi(b)Li_m(\zeta_N^b),$$
il faut remplacer $\phi_m^{CW}(N)(u).L_p(\chi^{-1}\omega^{1-m},m)$, dans la formule syntomique (4) précédente, par $\phi_m^{CW}(M)(u).L_p(\chi^{-1}\omega^{1-m},m)(N/M)^{m-1}.\mathrm{Eul}(m,\chi)$, où $\mathrm{Eul}(m,\chi)$ désigne le produit des facteurs euleriens en $l$ qui sont apparus plus haut.
\\--dans la formule étale (3) précédente, la relation de corestriction entre éléments de Deligne-Soulé ([HK],3.1.4.) montre qu'il faut remplacer $\chi_m^S(N)(u)$ par $\chi_m^S(N)(u).\mathrm{Eul}(m,\chi)$.
\\La commutativité du triangle des régulateurs montre alors que la formule (5) reste valable pour le corps $\mathbb{Q}(\zeta_M)^+$, quitte à multiplier $\mathcal{F}_{\psi}$ par un inversible, mais sans changer de notation.
\\Revenant au cas général, plongeons $F$ dans $\mathbb{Q}(\zeta_M)^+$, avec $M$ le conducteur de $F$. La corestriction de $\mathbb{Q}(\zeta_M)^+$ à $F$ se fait sans facteur parasite puisque les conducteurs sont les mêmes. La formule du théorème est ainsi démontrée pour $\epsilon\in\mathcal{U}_{\infty}^{\psi}$. Pour passer de $\epsilon\in\mathcal{U}_{\infty}^{\psi}$ à $\rho\in\mathfrak{X}_{\infty}^{\psi}$, on utilise le même argument de densité que dans la preuve du théorème 4.1.1: la restriction $\mathrm{Hom}_{G_{\infty}}(\mathfrak{X}_{\infty}^{\psi},\mathbb{Z}_p(m))\rightarrow \mathrm{Hom}_{G_{\infty}}(\mathcal{U}_{\infty}^{\psi},\mathbb{Z}_p(m))$ est injective pour presque tout $m$ ([KN], thm. 3.1.).
\end{proof}
Il reste maintenant à généraliser la mesure de Ihara-Kaneko-Yukinari (définition-lemme 4.2.1).

\begin{DefLemme}
Soit $\psi=\chi \omega^{i}$ comme précédemment, et $N=N_{\chi}$. Pour tout $\rho\in\mathfrak{X}^{\psi}_{\infty}$, pour tout $n\geq 1$, et tout $\bar{a}\in (\mathbb{Z}/p^n\mathbb{Z})^{*}$, fixons $b_{n,\bar{a}}\in\mathbb{Z}_p[\psi]$ tel que $\zeta_{p^n}^{b_{n,\bar{a}}}=(\sigma_a e_n^{1/p^n})^{ \rho-1}$, où $e_n$ est la $p$-unité $N_{\mathbb{Q}(\zeta_{f p^n})/K_n^{+}}(1-\zeta_{f p^n})$ qui intervient dans la définition des caractères de Soulé (voir les rappels ). Posons $c=a N^{-1}$ et 
$$\delta_n(\rho)=\displaystyle{\sum_{a=1,\ (a,p)=1}^{p^n-1}}b_{n,\bar{a}} c^{-1}\sigma_c \in \mathbb{Z}_p[\psi][G_{\infty}].$$
Alors $\delta_{\psi}(\rho)=\displaystyle{\lim_{n\rightarrow\infty}} \delta_n(\rho)\in\mathbb{Z}_p[\psi][[G_{\infty}]]$ existe et $$\delta_{\psi}(\rho).(1+T)=\displaystyle{\sum_{m\geq 1, m \equiv i (p-1)}}\frac{N^{m-1}\chi_m(\rho)}{m!}X^m, \ \mathrm{exp}X=(1+T).$$
\end{DefLemme}
\begin{proof}[\normalfont{\bfseries{Preuve.}}]
Ce sont les mêmes calculs que dans la définition-lemme 4.2.1.
\end{proof}

Dans ce qui suit, on considérera $\delta_{\psi}(\rho)$, pour $\rho\in\mathfrak{X}_{\infty}^{\psi}$, comme un élément de $\Lambda_{\psi}=\widehat{\mathcal{O}}_F[[G_{\infty}]]^{\psi}$.

\begin{Theoreme}
Pour tout caractère impair $\psi$ de $\mathrm{Gal}(K/\mathbb{Q})$, pour tout $\rho\in \mathfrak{X}_{\infty}^{\psi}$, on a: $\mathrm{Mel}(\mathcal{F}_{\psi}.\rho)=\delta_{\psi}(\rho)$ dans $\Lambda_{\psi}$. En particulier, $\widetilde{\mathcal{G}}_{\infty}^{\psi}$ est composé des $\delta_{\psi}(\rho)$, $\rho$ parcourant $\mathfrak{X}_{\infty}^{\psi}$. Si $\psi\neq \omega$, les $(\delta_{\psi}(\rho))^{\sharp}$ constituent l'idéal de Fitting sur $\Lambda_{\psi^*}$ du module $X_{\infty}^{'\psi^*}$.  
\end{Theoreme}
\begin{proof}[\normalfont{\bfseries{Preuve.}}]
Par définition de l'homomorphisme de Coleman et des caractères de Coates-Wiles, $$\mathrm{Col}(\epsilon)=\displaystyle{\sum_{m\geq 1}}(1-\chi(p)p^{m-1})\phi_m^{CW}(\epsilon)\frac{X^m}{m!}$$
pour tout $\epsilon\in\mathcal{U}_{\infty}^{\psi}$, et donc, d'après le Théorème 5.3.1. et la remarque 2 du § 5.1, $\mathrm{Mel}(\mathcal{F}_{\psi}.\rho)=\delta_{\psi}(\rho)$ pour tout $\rho\in\mathfrak{X}_{\infty}^{\psi}$.
\end{proof}

Pour avoir un résultat d'annulation galoisienne au niveau fini, introduisons quelques notations:
\\$\forall n\geq 1, \ \forall\rho\in\mathfrak{X}_{\infty}^{\psi},\ \forall\tau\in G_n$, soit $\bar{b}_{n,\tau}\in\mathbb{Z}/p^n\mathbb{Z}[\psi]$ tel que $\zeta_{p^n}^{\bar{b}_{n,\tau}}=(\tau(e_n)^{1/p^n})^{\rho-1}$, où $e_n$ est la $p$-unité spéciale qui intervient dans la définition des caractères de Soulé. Posons $$\overline{\mathfrak{d}}_n(\rho)=\displaystyle{\sum_{\tau\in G_n}}\bar{b}_{n,\tau} \tau^{-1}\in \mathbb{Z}_p[\psi^*]/p^n[G_n].$$

\begin{Corollaire}
Si $\psi(p)\neq 1$ et $\psi^*(p)\neq 1$, alors les $\sigma_N \overline{\mathfrak{d}}_n(\rho)$, $\rho$ parcourant $\mathfrak{X}_{\infty}^{\psi}$, composent l'idéal de Fitting de $A^{\psi^*}_n/p^n$ sur $\mathbb{Z}_p[\psi^*]/p^n[G_n]$.
\end{Corollaire}
\begin{proof}[\normalfont{\bfseries{Preuve.}}]
Dans les calculs justifiant l'existence de $\delta_{\psi}(\rho)$, on a en fait montré que $\delta_{\psi}(\rho)=\displaystyle{\lim_{\longleftarrow}}\pi_n(\delta_n(\rho))$, où $\pi_n$ est la projection de $\widehat{\mathcal{O}}_F[[G_{\infty}]]$ dans $\widehat{\mathcal{O}}_F/p^n[G_n]$. Or 
$$\pi_n(\delta_n(\rho)^{\sharp})=\displaystyle{\sum_{a\in(\mathbb{Z}/p^n\mathbb{Z})^{*}}} \bar{b}_{n,a} \sigma_{a^{-1} N}=\sigma_{N} \overline{\mathfrak{d}}_n(\rho)$$
(rappelons que $N=N_{\chi}$ est premier à $p$). On applique alors le corollaire 4 du théorème 3.4.2.
\end{proof}

\subsection{Calculs globaux}

Dans cette sous-section, on cherche à faire une étude "globale" (par opposition à l'étude "caractère par caractère" des sections précédentes) des parties "+ et -" de certains modules galoisiens ($p\neq 2$). \textit{L'hypothèse de semi-simplicité n'est plus nécessaire.} 
Toutes les notations antérieures ne faisant pas intervenir de caractère seront conservées ; les autres seront modifiées de façon naturelle. En particulier (cf. § 3.2) $\mathrm{tor}_{\Lambda} \mathcal{U}_{\infty}\simeq \mathbb{Z}_p[\Delta/\Delta_p](1)$ (cf. §5.1), $\widetilde{\mathcal{U}}_{\infty}=\mathcal{U}_{\infty}/\mathrm{tor}_{\Lambda} \mathcal{U}_{\infty}$, $W_{\infty}=\mathrm{tor}_{\Lambda} \mathcal{U}_{\infty}/\mathbb{Z}_p(1)$, $\widetilde{\mathfrak{X}}_{\infty}=\mathfrak{X}_{\infty}/W_{\infty}$ (c'est le \textbf{module de Bertrandias-Payan}, voir §6). Tous ces modules, ainsi que leurs parties plus et moins, sont des modules sur l'algèbre complète $\mathbb{A}=\mathbb{Z}_p[\Delta][[G_{\infty}]]$, qui contient la sous-algèbre d'Iwasawa $\Lambda=\mathbb{Z}_p[[\Gamma]]$ et qui peut être identifiée à $\Lambda[\mathrm{Gal}(K/\mathbb{Q})]$ (les groupes de Galois sont définis au début du §5). On notera $\widetilde{\mathrm{Mel}}:\ \widetilde{\mathcal{U}}_{\infty}\hookrightarrow\mathbb{A}$ l'homomorphisme injectif déduit de la suite exacte de Coleman (§5.1).

\begin{Theoreme}
Soit $\mathcal{F}$ la série caractéristique de $X_{\infty}^{-}$ et notons $[\mathcal{F}]:\ \widetilde{\mathfrak{X}}_{\infty}^{-}\rightarrow\mathbb{A}^{-}$ le morphisme composé $\ \widetilde{\mathrm{Mel}}\ \circ$ (multiplication par $\mathcal{F}$). On a une suite exacte de $\mathbb{A}$-modules:
$$0\rightarrow \alpha(X^{'}_{\infty}(-1))^{-}\rightarrow X_{\infty}^{-} \rightarrow \mathbb{A}^{-}/[\mathcal{F}]\widetilde{\mathcal{U}}_{\infty}^{-}\rightarrow \mathbb{A}^{-}/[\mathcal{F}]\widetilde{\mathfrak{X}}_{\infty}^{-}\rightarrow 0,$$	
\end{Theoreme}

\begin{proof}[\normalfont{\bfseries{Preuve.}}]
La même démonstration que dans le théo. 3.3.1. (mais sans passer par les caractères) donne une suite exacte de $\Lambda$-modules:
$$0\rightarrow \alpha(X^{'}_{\infty}(-1))^{-}\rightarrow X_{\infty}^{-} \rightarrow \widetilde{\mathcal{U}}_{\infty}^{-}/\mathcal{F}\widetilde{\mathcal{U}}_{\infty}^{-}\rightarrow \widetilde{\mathcal{U}}_{\infty}^{-}/\mathcal{F}\widetilde{\mathfrak{X}}_{\infty}^{-}\rightarrow 0.$$
Dans nos hypothèses abéliennes, la multiplication par $\mathcal{F}$ est $\mathbb{A}$-équivariante, donc le morphisme $\widetilde{\mathrm{Mel}}$ donne par passage au quotient un diagramme commutatif
\xymatrix{ &  & \widetilde{\mathcal{U}}_{\infty}^{-}/\mathcal{F} \widetilde{\mathcal{U}}_{\infty}^{-} \ar@{->>}[r] \ar@{^(->}[d]^{\widetilde{\mathrm{Mel}}} & \widetilde{\mathcal{U}}_{\infty}^{-}/\mathcal{F} \widetilde{\mathfrak{X}}_{\infty}^{-} \ar@{^(->}[d]^{\widetilde{\mathrm{Mel}}} 
\\ &  & \mathbb{A}^{-}/[\mathcal{F}] \widetilde{\mathcal{U}}_{\infty}^{-} \ar@{->>}[r] & \mathbb{A}^{-}/[\mathcal{F}] \widetilde{\mathfrak{X}}_{\infty}^{-}  
}
\ \\Les deux flèches verticales sont évidemment injectives et ont pour conoyau $\mathbb{A}^{-}/\widetilde{\mathrm{Mel}} \widetilde{\mathcal{U}}_{\infty}^{-}$. Le lemme du serpent donne alors un isomorphisme entre les noyaux des deux flèches horizontales, d'où la suite exacte du théorème.
\end{proof}

\begin{Corollaire}
Notons $M_{\infty}$ et $L_{\infty}$ les extensions de $K_{\infty}$ telles que $\mathfrak{X}_{\infty}=\mathrm{Gal}(M_{\infty}/K_{\infty})$ et $X_{\infty}=\mathrm{Gal}(L_{\infty}/K_{\infty})$. Soit $T_{\infty}$ le sous-corps de $M_{\infty}$ fixé par $\mathrm{tor}_{\Lambda} \mathfrak{X}_{\infty}$ (voir §6 ci-dessous). Alors $\mathrm{Gal}(L_{\infty}\cap T_{\infty}/K_{\infty})\simeq [\mathcal{F}]\widetilde{\mathfrak{X}}^-_{\infty}/[\mathcal{F}]\widetilde{\mathcal{U}}^-_{\infty}$ en tant que $\mathbb{A}$-modules.
\end{Corollaire} 

NB: d'après la conjecture de Greenberg, $\mathrm{Gal}(L_{\infty}\cap T_{\infty}/K_{\infty})$ devrait être égal à $X_{\infty}^{-}$ (voir §6).

\begin{proof}[\normalfont{\bfseries{Preuve.}}]
D'après l'appendice (§6 ci-dessous), $\mathrm{Gal}(L_{\infty}\cap T_{\infty}/K_{\infty})$ est naturellement isomorphe à $X_{\infty}^-/\alpha(X_{\infty}^{'}(-1))^-$. La suite exacte du théorème 5.4.1. permet de conclure.
\end{proof}

\begin{Theoreme}
Avec les hypothèses et notations du théorème 5.4.1., on a une relation entre $\mathbb{A}$-idéaux:
$$([\mathcal{F}]\widetilde{\mathfrak{X}}_{\infty}^{-})^{\sharp}=\mathrm{Fitt}_{\mathbb{A}}(\mathbb{A}\otimes_{\Lambda} X^{' +}_{\infty}).(\widetilde{\mathrm{Mel}}\widetilde{\mathcal{U}}_{\infty}^{-})^{\sharp}\subset \mathrm{Fitt}_{\mathbb{A}}(\mathbb{A}\otimes_{\Lambda}X^{'+}_{\infty}).$$
Il en résulte en particulier que $([\mathcal{F}]\widetilde{\mathfrak{X}}_{\infty}^{-})^{\sharp}\subset \mathrm{Ann}_{\mathbb{A}}(X^{'+}_{\infty}).$
\end{Theoreme}

\begin{proof}[\normalfont{\bfseries{Preuve.}}]
Nous allons faire la démonstration en plusieurs étapes:
\\1. La suite exacte de $\Lambda$-modules au début de la démonstration du théorème 5.4.1. donne une égalité de séries caractéristiques:	
$$\mathcal{F}.\mathrm{sc}(\widetilde{\mathcal{U}}_{\infty}^{-}/\mathcal{F}\widetilde{\mathfrak{X}}_{\infty}^{-})=\mathrm{sc}(\alpha(X^{'}_{\infty}(-1))^-).\mathrm{sc}(\widetilde{\mathcal{U}}_{\infty}^{-}/\mathcal{F}\widetilde{\mathcal{U}}_{\infty}^{-}).$$ 
Or, d'après le théorème de structure, $\widetilde{\mathcal{U}}_{\infty}^-$ s'injecte dans un $\Lambda$-module libre $\Lambda^r$ (où $r$ est le $\Lambda$-rang de $\widetilde{\mathcal{U}}_{\infty}^-$), avec un conoyau fini. Comme $p$ ne divise pas $\mathcal{F}$ (nullité de l'invariant $\mu$ de $X_{\infty}^-$), on en déduit facilement un pseudo-isomorphisme injectif $\widetilde{\mathcal{U}}_{\infty}^-/\mathcal{F} \widetilde{\mathcal{U}}_{\infty}^- \hookrightarrow \Lambda^r/\mathcal{F} \Lambda^r$, qui entraîne en particulier que le premier quotient n'a pas de sous-module fini non nul. On montre aussi, exactement comme dans le théorème 3.4.2. que les sous-modules finis maximaux de $\widetilde{\mathcal{U}}_{\infty}^-/\mathcal{F}\widetilde{\mathfrak{X}}_{\infty}^-$ et $X^{' +}_{\infty}$ sont entre eux en dualité de Kummer. En appliquant le lemme 3.4.1., on obtient alors une égalité de $\Lambda$-idéaux:
$$(\mathcal{F}^{\sharp}\Lambda).\mathrm{Fitt}_{\Lambda}(\widetilde{\mathcal{U}}_{\infty}^-/\mathcal{F}\widetilde{\mathfrak{X}}_{\infty}^-)^{\sharp}=\mathrm{Fitt}_{\Lambda}(X^{' +}_{\infty}).\mathrm{Fitt}_{\Lambda}(\widetilde{\mathcal{U}}_{\infty}^-/\mathcal{F}\widetilde{\mathcal{U}}_{\infty}^-)^{\sharp}.$$
On a montré en fait que le dernier idéal de Fitting à droite est $(\mathcal{F}^{\sharp})^r\Lambda$, mais on n'utilisera pas cette propriété. Il est préférable d'appliquer $\widetilde{\mathrm{Mel}}$ et d'utiliser la fonctorialité des idéaux de Fitting par rapport à l'extension des scalaires pour obtenir une égalité de $\mathbb{A}$-idéaux:
$$(\mathcal{F}^{\sharp}\mathbb{A}).\mathrm{Fitt}_{\mathbb{A}}(\widetilde{\mathrm{Mel}}\widetilde{\mathcal{U}}_{\infty}^-/[\mathcal{F}]\widetilde{\mathfrak{X}}_{\infty}^-)^{\sharp}=\mathrm{Fitt}_{\mathbb{A}}(\mathbb{A}\otimes_{\Lambda}X^{' +}_{\infty}).\mathrm{Fitt}_{\mathbb{A}}(\widetilde{\mathrm{Mel}}\widetilde{\mathcal{U}}_{\infty}^-/[\mathcal{F}]\widetilde{\mathcal{U}}_{\infty}^-)^{\sharp}.$$
2. Revenons à la suite exacte de Coleman (voir le début du §5.1):
$$0\rightarrow \widetilde{\mathcal{U}}_{\infty}^{-}\stackrel{\widetilde{\mathrm{Mel}}}{\rightarrow}\mathbb{A}^-\rightarrow V:=(\widehat{\mathcal{O}}_F/(\sigma_p-1)\widehat{\mathcal{O}}_F (1))^-\rightarrow 0$$
d'où nous déduisons une suite exacte:
$$0\rightarrow \widetilde{\mathcal{U}}_{\infty}^-/\mathcal{F} \widetilde{\mathcal{U}}_{\infty}^- \stackrel{\widetilde{\mathrm{Mel}}}{\rightarrow} \mathbb{A}^-/[\mathcal{F}]\widetilde{\mathcal{U}}_{\infty}^- \rightarrow V \rightarrow 0.$$
En tant que $\Lambda$-modules, on sait que $\mathbb{A}^-=\Lambda^r$ et $V$ (qui est un produit d'un certain nombre de copies de $\mathbb{Z}_p(1)$) est sans $\mathbb{Z}_p$-torsion, donc de $\Lambda$-dimension projective inférieure ou égale à un. Une propriété bien connue des idéaux de Fitting (voir e.g. [H], chap.3) donne alors une égalité de $\Lambda$-idéaux:
$$\mathrm{Fitt}_{\Lambda}(\Lambda^r/\mathcal{F}\widetilde{\mathcal{U}}_{\infty}^{-})=\mathrm{Fitt}_{\Lambda}(\widetilde{\mathcal{U}}_{\infty}^{-}/\mathcal{F}\widetilde{\mathcal{U}}_{\infty}^{-})
.\mathrm{Fitt}_{\Lambda}V.$$
En appliquant $\widetilde{\mathrm{Mel}}$ et par extension des scalaires, on obtient:
$$[\mathcal{F}]\widetilde{\mathcal{U}}_{\infty}^{-}=\mathrm{Fitt}_{\mathbb{A}}(\mathbb{A}^-/[\mathcal{F}]\widetilde{\mathcal{U}}_{\infty}^{-})=\mathrm{Fitt}_{\mathbb{A}}(\widetilde{\mathrm{Mel}} \widetilde{\mathcal{U}}_{\infty}^{-}/[\mathcal{F}]\widetilde{\mathcal{U}}_{\infty}^{-})
.\mathrm{Fitt}_{\mathbb{A}}(\mathbb{A}\otimes_{\Lambda}V).$$
Or $\widetilde{\mathrm{Mel}}\widetilde{\mathcal{U}}_{\infty}^-=\mathrm{Fitt}_{\mathbb{A}}V$ d'après la suite de Coleman, d'où $$[\mathcal{F}]\widetilde{\mathcal{U}}_{\infty}^{-}=\mathrm{Fitt}_{\mathbb{A}}(\widetilde{\mathrm{Mel}}\widetilde{\mathcal{U}}_{\infty}^-/[\mathcal{F}]\widetilde{\mathcal{U}}_{\infty}^-).\widetilde{\mathrm{Mel}}\widetilde{\mathcal{U}}_{\infty}^-.$$
Le même raisonnement exactement, en remplaçant $\mathcal{F}\widetilde{\mathcal{U}}_{\infty}^{-}$ par $\mathcal{F}\widetilde{\mathfrak{X}}_{\infty}^{-}$, donne:
$[\mathcal{F}]\widetilde{\mathfrak{X}}_{\infty}^{-}=\mathrm{Fitt}_{\mathbb{A}}(\widetilde{\mathrm{Mel}}\widetilde{\mathcal{U}}_{\infty}^-/[\mathcal{F}]\widetilde{\mathfrak{X}}_{\infty}^-).\mathrm{Fitt}_{\mathbb{A}}(\mathbb{A}\otimes_{\Lambda}V)$.
\item En multipliant la relation de 1. par $\mathrm{Fitt}_{\mathbb{A}}(\mathbb{A}\otimes_{\Lambda}V)$ et en tenant compte des deux dernières relations de 2., on obtient une égalité de $\mathbb{A}$-idéaux:
$$(\mathcal{F}^{\sharp} \mathbb{A}).([\mathcal{F}]\widetilde{\mathfrak{X}}_{\infty}^-)^{\sharp}=\mathrm{Fitt}_{\mathbb{A}}(\mathbb{A}\otimes_{\Lambda}X^{' +}_{\infty}).([\mathcal{F}]\widetilde{\mathcal{U}}_{\infty}^{-})^{\sharp}.$$
Or $[\mathcal{F}]\widetilde{\mathcal{U}}_{\infty}^{-}=\widetilde{\mathrm{Mel}}(\mathcal{F}\widetilde{\mathcal{U}}_{\infty}^-)=(\mathcal{F}\mathbb{A}).\widetilde{\mathrm{Mel}}\widetilde{\mathcal{U}}_{\infty}^{-}$, d'où
$$(\mathcal{F}^{\sharp} \mathbb{A}).([\mathcal{F}]\widetilde{\mathfrak{X}}_{\infty}^-)^{\sharp}=(\mathcal{F}^{\sharp}\mathbb{A}).\mathrm{Fitt}_{\mathbb{A}}(\mathbb{A}\otimes_{\Lambda}X^{' +}_{\infty}).(\widetilde{\mathrm{Mel}}\widetilde{\mathcal{U}}_{\infty}^{-})^{\sharp}.$$
$\mathcal{F}^{\sharp}$ est un élément de $\Lambda$, donc n'est pas un diviseur de zéro dans $\mathbb{A}$ ; on en déduit immédiatement la première relation énoncée dans le théorème. Il en résulte en particulier que pour tout $\delta\in([\mathcal{F}]\widetilde{\mathfrak{X}}_{\infty}^{-})^{\sharp}$, on a $\delta(1\otimes_{\Lambda}X_{\infty}^{'+})=\delta\otimes_{\Lambda}X_{\infty}^{'+}=(0)$. La propriété universelle du produit tensoriel entraîne que $\delta X^{'+}_{\infty}=(0)$, i.e. $\delta\in\mathrm{Ann}_{\mathbb{A}}X_{\infty}^{'+}$.

\end{proof}

Pour compléter le théorème 5.4.1. il faudrait donner une description explicite de l'idéal $[\mathcal{F}]\widetilde{\mathfrak{X}}_{\infty}^-$ dans le style des sous-sections précédentes du §5. Une telle description "globale" (par opposition à "caractère par caractère") nécéssiterait probablement d'introduire des techniques de la théorie d'Iwasawa équivariante. Nous espérons y revenir dans un travail ultérieur.

\subsection{Comparaison avec des résultats de Solomon}

Dans son exposé [Sol1] et sa prépublication ultérieure [Sol2], D.Solomon a construit certains éléments annulateurs de la partie "plus" du groupe de classes en procédant par montée (alors que nous avons procédé par descente). Résumons sa démarche:
\\1. Sous les hypothèses du §5.4, Solomon construit pour tout $n\geq 1$ un homomorphisme $$\overline{\mathfrak{d}}_n:\ (\mathfrak{X}_n^{-}/p^n)^{\sharp}\rightarrow\mathbb{Z}/p^n\mathbb{Z}[G_n]^+$$
	 qui (après traduction) est l'analogue (sans les caractères) du $\overline{\mathfrak{d}}_n$ du corollaire du théorème 5.3.2.
\\2. Par la méthode de Thaine, il montre que $\overline{\mathfrak{d}}_n(\rho)$ annule $A_n^+/p^n$, puis, par passage à la limite projective le long de la tour cyclotomique, il obtient un homomorphisme $\mathfrak{d}_{\infty}:\ (\mathfrak{X}_{\infty}^{-})^{\sharp}\rightarrow\mathbb{A}^+$ dont l'image annule $X^+_{\infty}$, ainsi qu'une suite exacte analogue (après traduction) à celle du théorème 5.4.1.
\\3. Dans le cas particulier de $\mathbb{Q}(\zeta_p)$, en redescendant $\mathrm{Im}\mathfrak{d}_{\infty}$ au niveau $n$, il détermine l'idéal de Fitting du dual de Pontryagin du quotient des unités par
 les unités cyclotomiques. Modulo un résultat de [CG], c'est équivalent à notre proposition 4.2.1.	  

\section{Appendice sur le "miroir" (Spiegelung)}

Dans cet appendice, on ne fait \textit{aucune hypothèse sur K autre que} $\mu_p\subset K$.
Les relations du "miroir" étant une combinaison d'isomorphisme (corps de classes) et de dualité (Kummer), c'est au-dessus de $K_{\infty}$, dans le cadre de la théorie d'Iwasawa, qu'elles s'expriment le mieux. Voici le schéma galoisien qui contient tous les résultats de Spiegelung ([Iw2], [N2], [LMN]):
\\
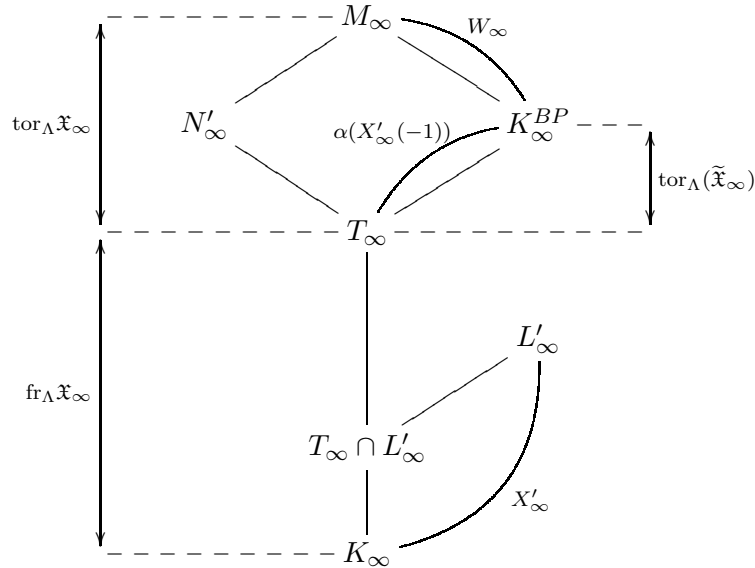
\xymatrix{ 
& & \ar@{--}[rr]& & M_{\infty} \ar@{-}[dr] \ar@{-}[dl] \ar@/^1pc/@{-}[dr]^{W_{\infty}} & &  \\ 
& & & N'_{\infty} \ar@{-}[dr] & & K_{\infty}^{BP} \ar@{-}[dl] \ar@{--}[r] \ar@/_1pc/@{-}[dl]_{\alpha(X'_{\infty}(-1))}&  \ar@{<->}[d]^{\mathrm{tor}_{\Lambda}(\widetilde{\mathfrak{X}}_{\infty})} \\
& & \ar@{<->}[uu]^{\mathrm{tor}_{\Lambda} \mathfrak{X}_{\infty}} \ar@{<->}[ddd]_{\mathrm{fr}_{\Lambda} \mathfrak{X}_{\infty}} \ar@{--}[rr] & & T_{\infty} \ar@{-}[dd] \ar@{--}[rr]& & \\
& & & & & L'_{\infty} \ar@{-}[dl] & \\
& & & & T_{\infty}\cap L'_{\infty} \ar@{-}[d] & & \\
& & \ar@{--}[rr] & & K_{\infty} \ar@/_2pc/@{-}[uur]_{X'_{\infty}} & & }     
\ \\
Les notations pour la plupart des extensions sont celles de [Iw 2]: $M_{\infty}$ est la pro-$p$-extension $p$-ramifiée abélienne maximale de $K_{\infty}$, $N'_{\infty}$ est l'extension de $K_{\infty}$ obtenue en ajoutant toutes les racines $p^k$-ièmes de toutes les $p$-unités de $K_{\infty}$, $L'_{\infty}$ est la pro-$p$-extension abélienne non ramifiée $p$-décomposée maximale de $K_{\infty}$; $\mathfrak{X}_{\infty}=\mathrm{Gal}(M_{\infty}/K_{\infty})$, $X'_{\infty}=\mathrm{Gal}(L'_{\infty}/K_{\infty})$; $T_{\infty}$ est le corps fixe de $\mathrm{tor}_{\Lambda}\mathfrak{X}_{\infty}$, donc $\mathrm{Gal}(T_{\infty}/K_{\infty})\simeq \mathrm{fr}_{\Lambda} \mathfrak{X}_{\infty}:=\mathfrak{X}_{\infty}/ \mathrm{tor}_{\Lambda} \mathfrak{X}_{\infty}$. Notons que la conjecture faible de Leopoldt (valable pour $K_{\infty}$) entraîne $\mathfrak{X}_{\infty}^{+}=(\mathrm{tor}_{\Lambda}\mathfrak{X}_{\infty})^{+}$ et $\mathrm{fr}_{\Lambda}\mathfrak{X}_{\infty}=(\mathrm{fr}_{\Lambda} \mathfrak{X}_{\infty})^{-}$. En utilisant la théorie de Kummer, Iwasawa montre dans [Iw2] que:
\\1. $\mathrm{Gal}(M_{\infty}/N'_{\infty})\simeq\alpha(X'_{\infty}(-1))$, où $\alpha(.)$ désigne l'adjoint.
\\2. $\mathrm{Gal}(N'_{\infty}/T_{\infty})\approx W_{\infty}:=\oplus_{v|p}(\mathrm{Ind}_{\Lambda_v}^{\Lambda} \mathbb{Z}_p(1))/\mathbb{Z}_p(1)$, où $v$ parcourt les $p$-places de $K_{\infty}$ et $\Lambda_v$ est l'algèbre d'Iwasawa locale en $v$ et $\approx$ signifie "pseudo-isomorphisme".
\\3. $\mathrm{fr}_{\Lambda} \mathfrak{X}_{\infty} = \mathrm{fr}_{\Lambda}\mathrm{Gal}(N'_{\infty}/K_{\infty})$.

Ces résultats peuvent être affinés en introduisant, comme dans [N1], [N2], le \textbf{corps de Bertrandias-Payan} $K_{\infty}^{BP}$, dont on ne rappellera pas la définition mais qui peut être caractérisé comme étant le corps fixe de $W_{\infty}$. Alors, avec les notations du §3.2. , $\mathrm{Gal}(K_{\infty}^{BP}/K_{\infty})\simeq \widetilde{\mathfrak{X}}_{\infty}$ et l'on a (voir [LMN]):
\\1. $\mathrm{Gal}(M_{\infty}/N'_{\infty})\simeq\mathrm{Gal}(K_{\infty}^{BP}/K_{\infty}^{BP}\cap N'_{\infty})$ et $\mathrm{Gal}(K_{\infty}^{BP}/T_{\infty})\simeq\alpha(X'_{\infty}(-1)) \simeq \mathrm{tor}_{\Lambda}\widetilde{\mathfrak{X}}_{\infty}$
\\2. $\mathrm{Gal}(N'_{\infty}/N'_{\infty}\cap K_{\infty}^{BP})\simeq \mathrm{Gal}(M_{\infty}/K_{\infty}^{BP})=W_{\infty}$
\\3. $\mathrm{Gal}(N'_{\infty}\cap K_{\infty}^{BP}/T_{\infty})$ est fini.

Autrement dit, $\mathrm{tor}_{\Lambda} \mathfrak{X}_{\infty}$ contient un sous-module d'indice fini isomorphe à $W_{\infty}\oplus \alpha(X'_{\infty}(-1))$. Ces résultats proviennent essentiellement d'une suite exacte de Poitou-Tate
$$0\rightarrow W_{\infty} \rightarrow \mathrm{tor}_{\Lambda} \mathfrak{X}_{\infty} \rightarrow \mathrm{tor}_{\Lambda} \widetilde{\mathfrak{X}}_{\infty}\simeq \alpha(X'_{\infty}(-1)) \rightarrow 0$$ 
qui exprime le fait que la surjection canonique $\mathfrak{X}_{\infty} \rightarrow X'_{\infty}$ induit un homomorphisme $\mathrm{tor}_{\Lambda} \mathfrak{X}_{\infty} \rightarrow X'_{\infty}$ qui se factorise à travers $\eta_{\infty}:\ \mathrm{Gal}(K_{\infty}^{BP}/T_{\infty})=\mathrm{tor}_{\Lambda} \widetilde{\mathfrak{X}}_{\infty} \rightarrow X'_{\infty}$. Alors (voir [N2]):
\\1. $\mathrm{Gal}(M_{\infty}/N'_{\infty})\simeq\mathrm{Gal}(K_{\infty}^{BP}/K_{\infty}^{BP}\cap N^{'}_{\infty})$
\\2. $\mathrm{Gal}(K_{\infty}^{BP}/T_{\infty})\simeq \alpha(X^{'}_{\infty}(-1))\simeq \mathrm{tor}_{\Lambda} \widetilde{\mathfrak{X}}_{\infty}$  
\\ Supposons en outre que $K$ est CM:
\\3. $\mathrm{Ker}\eta_{\infty}^{-}=(0)$ et $(\mathrm{tor}_{\Lambda} \widetilde{\mathfrak{X}}_{\infty})^{-}\simeq \alpha(X'_{\infty}(-1))^{-}$
\\4. $\mathrm{Im}\eta_{\infty}=\mathrm{Gal}(L'_{\infty}/L'_{\infty}\cap T_{\infty})\simeq X_{\infty}^{'+}\oplus \alpha(X'_{\infty}(-1))^{-}$  
\\(certains résultats de [N2] sont exprimés à pseudo-isomorphisme près, mais ils sont en fait valables à isomorphisme près, avec exactement les mêmes démonstrations). D'après les propriétés fonctorielles de l'adjoint, en notant $(.)^{\sharp}$ l'involution d'Iwasawa, $\alpha(X'_{\infty}(-1))^{-} \approx (X_{\infty}^{'+})^{\sharp}$. La conjecture de Greenberg prédit la finitude de $X_{\infty}^{'+}$, donc, puisqu'un adjoint n'a pas de sous-module fini non nul, la nullité de $\alpha(X'_{\infty}(-1))^{-}$. Par conséquent, elle équivaut à la finitude de $\mathrm{Im} \eta_{\infty} =\mathrm{Gal} (L'_{\infty}/L'_{\infty}\cap T_{\infty})$ et l'on peut la généraliser au cas non CM:
\begin{Conjecture}
Pour un corps de base $K$ quelconque, $[L'_{\infty}:L'_{\infty}\cap T_{\infty}]< \infty$
\end{Conjecture}

\textbf{Cas particulier du corps de base $\mathbb{Q}(\zeta_p)$:}
Dans ce cas, il n'y a qu'une seule $p$-place au-dessus de $p$ et elle est totalement ramifiée dans $\mathbb{Q}(\zeta_{p^{\infty}})/\mathbb{Q}$. On voit alors facilement que:
\\$M_{\infty}=K_{\infty}^{BP}$ et $\mathrm{tor}_{\Lambda}\mathfrak{X}_{\infty}=\alpha(X'_{\infty}(-1))$; $T_{\infty}=N'_{\infty}=N_{\infty}=$ l'extension de $\mathbb{Q}(\zeta_{p^{\infty}})$ obtenue en ajoutant toutes les racines $p^k$-ièmes de toutes les unités de $\mathbb{Q}(\zeta_{p^{\infty}})$; $L'_{\infty}=L_{\infty}=$ la pro-$p$-extension abélienne non ramifiée maximale de $\mathbb{Q}(\zeta_{p^{\infty}})$. La conjecture de Greenberg équivaut à $[L_{\infty}:L_{\infty}\cap T_{\infty}]<\infty$ et la conjecture de Vandiver à $[L_{\infty}:L_{\infty}\cap T_{\infty}]=1$.
\\ Pour tout $n\geq 0$, soit $C_n$ le groupe des unités cyclotomiques de $\mathbb{Q}(\zeta_{p^{n+1}})$, et soit $C_{\infty}=\displaystyle{\lim_{\longrightarrow}} C_n$. On sait que $(\mathcal{O}_n^{\times\ +}:C_n^{+})$ est égal au nombre de classes de $\mathbb{Q}(\zeta_{p^{n+1}})^{+}$, et la finitude de ces indices entraîne immédiatement que $T_{\infty}=N_{\infty}$ coïncide avec l'extension $Cyclo$ obtenue en ajoutant à $\mathbb{Q}(\zeta_{p^{\infty}})$ toutes les racines $p^k$-ièmes de tous les éléments de $C_{\infty}$. Or, compte tenu de la définition des caractères de Soulé et du théorème 4.1.1. de la partie 4, il n'est pas difficile de montrer que $Cyclo$ coïncide avec le sous-corps de $M_{\infty}$ fixé par le noyau de la représentation d'Ihara $\mathrm{Ih}:\ \mathfrak{X}_{\infty}\rightarrow \mathcal{A}^{\times}$ (pour les détails, voir [IS], prop.5). Notant $M_{\infty}^{-}$ le sous-corps fixé par $\mathfrak{X}_{\infty}^{+}$, on a:

\begin{Proposition}
$Cyclo \subset M_{\infty}^{-}$, et  $Cyclo=M_{\infty}^{-}$ si et seulement si la conjecture de Greenberg est vraie pour $\mathbb{Q}(\zeta_{p^{\infty}})$.
\end{Proposition} 
\begin{proof}[\normalfont{\bfseries{Preuve.}}]
Il est évident que $Cyclo \subset M_{\infty}^{-}$ puisque $(\mathrm{tor}_{\Lambda}\mathfrak{X}_{\infty})^{+}=\mathfrak{X}_{\infty}^{+}$. Or, d'après ce qui précède, la conjecture de Greenberg équivaut à $(\mathrm{tor}_{\Lambda}\mathfrak{X}_{\infty})^{-}=(0)$, i.e. à $\mathfrak{X}_{\infty}^{+}=\mathrm{tor}_{\Lambda}\mathfrak{X}_{\infty}$, ou encore à $M_{\infty}^{-}=T_{\infty}$. Mais ici $T_{\infty}=Cyclo$.
\end{proof}
\textbf{Remarque:} Ce résultat est à comparer avec la prop. 5 de [IS], qui utilise la conjecture de Vandiver.
\\ \ \\
\\ \ \\
Thong NGUYEN QUANG DO, Vésale NICOLAS
\\Université de Franche-Comté
\\Laboratoire de Mathématiques, CNRS UMR 6623
\\16 route de Gray, 25030 Besançon Cedex, France.
\\
\\thong.nguyen-quang-do@univ-fcomte.fr, vnicolas@univ-fcomte.fr.
\\ \ \\

\end{document}